\documentclass[a4paper,twoside]{amsart}
\usepackage{amsmath}
\usepackage{amssymb}
\usepackage{amsthm}
\usepackage{hyperref}
\usepackage{geometry}
\usepackage{graphicx}
\usepackage{subfig}
\usepackage{xcolor}

\begin{document}
\title{On the Stability of the solitary waves to the (Generalized) Kawahara Equation}
\author[A. Kabakouala and L. Molinet]{Andr\'e  Kabakouala and Luc Molinet}

\address{Andr\'e  Kabakouala, Laboratoire de Math\'ematiques et Physique Th\'eorique (CNRS UMR 7350), 
F\'ed\'eration Denis Poisson (FR CNRS 2964), Universit\'e Fran\c cois Rabelais-Tours, Parc Grandmont, 37200 Tours, France.}
\email{Andre.Kabakouala@lmpt.univ-tours.fr}
\address{Luc Molinet, Laboratoire de Math\'ematiques et Physique Th\'eorique (CNRS UMR 7350), 
F\'ed\'eration Denis Poisson (FR CNRS 2964), Universit\'e Fran\c cois Rabelais-Tours, Parc Grandmont, 37200 Tours, France.}
\email{Luc.Molinet@lmpt.univ-tours.fr}

% ENVIRONNEMENT THEOREME
\theoremstyle{plain}
\newtheorem{Theo}{Theorem}[section]
\newtheorem{Pro}{Proposition}[section]
\newtheorem{Lem}{Lemma}[section]
\newtheorem{Cor}{Corollary}[section]
\newtheorem{Def}{Definition}[section]
\newtheorem{Rem}{Remark}[section]
\newtheorem{Ex}{Example}[section]
\newtheorem*{Ack}{Acknowledgements}
\numberwithin{equation}{section}
\def\R{{\mathbb R}}
\def\N{{\mathbb N}}

% TITRE 
\maketitle

% ABSTRACT 
\begin{abstract}
In this paper we investigate the orbital stability of solitary waves to the (generalized) Kawahara equation (gKW) 
which is  a  fifth order dispersive equation.  For some values of the power of the nonlinearity, we prove the orbital stability 
in the energy space $H^{2}(\mathbb{R})$ of two branches of even solitary waves of gKW  by combining the well-known spectral
method introduced by Benjamin \cite{MR0338584}   
with continuity arguments. We construct the first family of even solitons 
by applying the implicit function theorem in the neighborhood of the explicit
solitons of gKW found by Dey et al. \cite{MR1427885}. The second family consists of even travelling waves with low speeds. 
They are solutions of a constraint minimization
problem on the line and rescaling of perturbations
of the soliton of gKdV with speed 1.
\end{abstract}
\vspace{1cm}

% INTRODUCTION 
\section{Introduction}

% EQUATION gKW
We consider the generalized Kawahara equation (gKW):
\begin{equation}\label{Int1}
\partial_{t}u+u^{p}\partial_{x}u+\partial^{3}_{x}u-\mu\partial^{5}_{x}u=0,
~~(t,x)\in\mathbb{R}^{*}_{+}\times\mathbb{R},
\end{equation}
with initial data $u(0)=u_{0}\in H^{2}(\mathbb{R})$,
where $p\in\mathbb{N}^{*}$ denotes the power of nonlinearity,
and $\mu>0$ the parameter which control the fifth-ordre dispersion term. 
In the case $ p=1$, this equation has been derived by Kawahara \cite{Kawa} as a model for
water waves with weak amplitude  in the long-wave regime approximation for moderate values of surface tension and a Weber number 
close to $1/3$. 
 For such
Weber numbers the usual description of long water waves via the Korteweg-de Vries (KdV)
equation fails since the cubic term in the linear dispersion relation vanishes and fifth order
dispersion becomes relevant at leading order.
Note that positive values of the parameter 
 $ \mu $ in \eqref{Int1} 
 correspond to Weber numbers larger  than $1/3$.

% PROBLEME DE CAUCHY 
The Cauchy problem associated to \eqref{Int1} is locally 
well-posed in $H^{2}(\mathbb{R})$ (see for instance  \cite{MR1044731}).
% LOIS DE CONSERVATION DE gKW
The $H^{2}$-solutions of \eqref{Int1} satisfy the following
two conservation laws in time:
\begin{equation}\label{Int01}
E_{p,\mu}(u(t))=\int_{\mathbb{R}}\left[\frac{\mu}{2}
(\partial^{2}_{x}u)^{2}(t)
+\frac{1}{2}(\partial_{x}u)^{2}(t)
-\frac{1}{(p+1)(p+2)}u^{p+2}(t)\right]=E_{\mu}(u_{0})
~~\text{(energy)}
\end{equation}
and
\begin{equation}\label{Int2}
V(u(t))=\frac{1}{2}\int_{\mathbb{R}}u^{2}(t)=V(u_{0})
~~\text{(mass)}.
\end{equation}
These conserved quantities enable to extend the solutions for all positive times so that \eqref{Int1} is actually globally well-posed in $ H^2(\R)$.

% EQUATION DES SOLITONS DE gKW
We are interested in  solitary waves of gKW, i.e., 
the solutions to equation \eqref{Int1} of the form
$u(t,x)=\varphi_{c,p,\mu}(x-ct)$,
traveling to the right with the speed $c>0$. 
Substituting $u$ by $\varphi_{c,p,\mu}$ in \eqref{Int1},
integrating on $ \R $  with the assumption
$\partial^{k}_{x}\varphi_{c,p,\mu}(\pm\infty)=0$ for $k=0,\ldots 4$, 
we obtain
\begin{equation}\label{Int3}
\mu\partial^{4}_{x}\varphi_{c,p,\mu}(x)
-\partial^{2}_{x}\varphi_{c,p,\mu}(x)
+c\varphi_{c,p,\mu}(x)=\frac{1}{p+1}\varphi^{p+1}_{c,p,\mu}(x),
~~\forall x\in\mathbb{R}.
\end{equation}

% SOLITONS EXPLICITES DE gKW
In  \cite{MR1427885}, Dey et al. compute 
 explicit solutions to \eqref{Int3}  that write : 
\begin{equation}\label{Int6}
\varphi_{c,p,\mu}(x)=\left[\frac{(p+1)(p+4)(3p+4)c}{8(p+2)}\right]^{1/p}
\text{sech}^{4/p}\left[\frac{p\sqrt{(p^{2}+4p+8)c}}{4(p+2)}x\right],
~~\forall x\in\mathbb{R},
\end{equation}
with
\begin{equation}\label{Int7}
c=\frac{2^{2}(p+2)^{2}}{(p^{2}+4p+8)^{2}\mu}.
\end{equation}
and study their  orbital stability. More generally,  the existence of solutions to \eqref{Int1} with $ \mu>0 $ and  has been proven in \cite{MR1667522} by solving  an associated constrained minimization problem. The orbital stability of the set of ground states to   this minimization problem is also 
 studied. Note that  \eqref{Int7} is a fourth order differential equation and, up to our knowledge, no  uniqueness result (up to symmetries) is known even for ground state solutions.
%Note that, others explicit solitons 
%can be found in the literature (for some nonlinearity $p$),
%but they will not be considered in this work 
%(see for instance Natali \cite{MR2602415}).
% SOLITONS DE gKdV

On the other hand, for  $\mu=0$,  \eqref{Int0001} becomes the well-known generalized Korteweg-de Vries
equation (gKdV):
\begin{equation}\label{Int0001}
\partial_{t}u+u^{p}\partial_{x}u+\partial^{3}_{x}u=0,
~~(t,x)\in\mathbb{R}^{*}_{+}\times\mathbb{R}
\end{equation}
and thus, for $ \mu=0 $  in \eqref{Int3},  we  recover the equation
of solitons of gKdV: 
\begin{equation}\label{Int0002}
-\partial^{2}_{x}\varphi_{c,p,0}(x)+c\varphi_{c,p,0}(x)
=\frac{1}{p+1}\varphi^{p+1}_{c,p,0}(x),
~~\forall x\in\mathbb{R}.
\end{equation}
Recall that, the solitons $\varphi_{c,p,0}$ of gKdV are unique,
up to  translations, and to the transformation: 
$\varphi_{c,p,0}\mapsto-\varphi_{c,p,0}$
if $p$ is even. Moreover, they are explicitly defined by:
\begin{equation}\label{Int0003}
\varphi_{c,p,0}(x)=\left[\frac{(p+1)(p+2)c}{2}\right]^{1/p}
\text{sech}^{2/p}\left[\frac{p\sqrt{c}}{2}x\right],
~~\forall x\in\mathbb{R}.
\end{equation}

% EXISTENCE 1: THEOREME DES FONCTIONS IMPLICITES 
In this  paper, we start by  constructing two branches of solutions to equation \eqref{Int3} with $ \mu=1$ and establishing some uniqueness results on these solutions.
Firstly, we fix $p\in\mathbb{N}^{*}$, $\mu=1$ 
and $c_{p}=\frac{2^{2}(p+2)^{2}}{(p^{2}+4p+8)^{2}}$.
By applying the Implicit Function Theorem 
in the neighborhood of the explicit solution $\varphi_{c_{p},p}$,
we construct a continuous in $ H^4(\R) $ branch $\{\varphi_{c,p},
~c\in]c_{p}-\delta_{p},c_{0}+\delta_{p}[\}$, with $0<\delta_{p}\ll 1$,
of even solutions to \eqref{Int3}. For each $ c\in ]c_{0}-\delta_{p},c_{0}+\delta_{p}[$, $\varphi_{c,p} $ is the unique even $ H^4(\R) $ solution of \eqref{Int3} in some 
 $ H^4$-neighborhood of $ \varphi_{c_p,p} $.
% EXISTENCE 2:  CONCENTRATION COMPACITE 
Secondly, for all $p\in\mathbb{N}^{*}$, $c>0$
and $\mu>0$, following \cite{MR1667522}, we  minimize in the even functions of $ H^2(\R) $ 
the functional: 
\begin{equation}\label{Int00001}
I_{c,\mu}(\psi)=\int_{\mathbb{R}}\left[\frac{\mu}{2}
(\partial^{2}_{x}\psi)^{2}
+\frac{1}{2}(\partial_{x}\psi)^{2}+\frac{c}{2}\psi^{2}\right],
\end{equation}
under the constraint:
\begin{equation}\label{Int00002}
K_{p}(\psi)=\frac{1}{(p+1)(p+2)}\int_{\mathbb{R}}\psi^{p+2}
=K_{p}(\varphi_{c,p,0}).
\end{equation}
For $ \mu>0 $ small enough we prove the uniqueness of the associated  even ground states by using ideas of  \cite{MR2859931} .
In this way, we construct a family $\{\varphi_{c,p,\mu}, 0<\mu\ll 1\}$ 
of even solutions to \eqref{Int3} such that 
\begin{equation}\label{Int00003}
\lim_{\mu\to 0^{+}}\left\|\varphi_{c,p,\mu}
-\varphi_{c,p,0}\right\|_{H^{1}(\mathbb{R})}=0,
\end{equation}
where $\varphi_{c,p,0}$ is defined by \eqref{Int0002}-\eqref{Int0003}.
 Noticing that $ u$ is a solution to \eqref{Int3} with $ c=1$, $ p=p_0\ge 1$ and  $\mu>0 $ if and only if $v=\mu^{1/p}u(\sqrt{\mu}\cdot) $ is a solution to 
 \eqref{Int3} with $c=\mu $, $p=p_0 $ and $ \mu=1 $,  we obtain a $ H^1$-continuous branch
 $$
 \{\varphi_{c,p,1},~c\in]0,\delta'_{p}[\} \quad\mbox{with}\quad 0<\delta'_{p}\ll 1,
 $$
of even solutions to \eqref{Int3} (with $\mu=1$) traveling with  low speeds. For each $ c\in  ]0,\delta'_{p}[$,  $ \varphi_{c,p} $ is the unique even solution   of the constraint minimizing problem \eqref{qq}.

% DEFINITION DE LA STABILITE ORBITALE
The main result of this paper is the orbital stability of the solitary waves that form these two branches. 
Note that this improves earlier results (see \cite{MR1667522}) where the stability of the set of ground states is proven. Before stating our main result, let us recall the definition of orbital stability and define what we will call ground state solutions to \eqref{Int3}.
\begin{Def}[Orbital Stability]\label{Def-Stab}
Let $\phi\in H^{2}(\mathbb{R})$ be a solution of \eqref{Int3}. 
We say that $\phi$ is orbitally stable in $H^{2}(\mathbb{R})$, 
if for all $\varepsilon>0$, there exists $\delta_{\varepsilon}>0$, 
such that for all initial data 
$u_{0}\in H^{2}(\mathbb{R})$, satisfying 
\begin{equation}\label{Int4}
\left\|u_{0}-\phi\right\|_{H^{2}(\mathbb{R})}\le\delta_{\varepsilon},
\end{equation}
the solution $u\in C(\mathbb{R}_{+},H^{2}(\mathbb{R}))$ of gKW 
emanating form $u_{0}$ satisfies
\begin{equation}\label{Int5}
\sup_{t\in\mathbb{R}_{+}}
\inf_{z\in\mathbb{R}}
\left\|u(t,\cdot+z)-\phi\right\|_{H^{2}(\mathbb{R})}\le\varepsilon.
\end{equation}
\end{Def}
In this paper, for $ s\ge 0 $, we set $H^s_{e}(\R)=\left\{u\in H^{s}(\mathbb{R}):u(-\cdot)=u(\cdot)\right\}$ with  
$L^2_{e}(\R)=H^0_e(\R) $. This space, endowed with  the metric of  
 $H^{s}(\mathbb{R})$ is an  Hilbert space.
\begin{Def}[Even ground state solution]\label{groundstate}
We say that a solution to \eqref{Int3} is an even ground state solution to \eqref{Int3},  if it is also a solution to the constraint minimization problem
\begin{equation}\label{qq}
S^{\beta}_{c,p,\mu}=\inf\left\{I_{c,\mu}(\psi):~\psi\in H^2_e(\R),
~K_{p}(\psi)=\beta\right\}.
\end{equation}
for some $ \beta>0$.
\end{Def}

% RESULTAT PRINCIPAL 
Our main results can be summarized as follows :  
\begin{Theo}\label{Main-Result}
Let us  fix $\mu=1$ and  set $c_p=\frac{2^{2}(p+2)^{2}}
{(p^{2}+4p+8)^{2}}$ for any $ p\ge 1$.\\
(i) For $p\in\left\{1,2,3,4\right\}$, there exist $\delta_{p}>0$ such that for any $ c\in ]c_p-\delta_p, c_p+\delta_p[ $,  \eqref{Int3} has a unique even solution  in the ball of $ H^4(\R) $ centered in $\varphi_{c_p,p,1}$  (explicitly defined in \eqref{Int6}-\eqref{Int7}) with radius $ \delta_p$.

 These solutions form  a curve of class $ C^1 $ in $H^{4}(\R)$, passing by $\varphi_{c_p,p,1}$, 
%\begin{equation}\label{Int8}
%\left\{\varphi_{c,p},~c\in]c(p)-\delta_{p},c(p)+\delta_{p}[\right\}
%\subset H^{4}(\mathbb{R}),
%\end{equation}
of even solitary waves to \eqref{Int1} which are all  stable. 
\\
(ii) For $p\in\left\{1,2,3\right\}$, there exist $\delta'_{p}>0$ such that, for any $ c\in ]0,\delta_p'[$, there exists a unique even ground state 
$\varphi_{c,p,1}$ to \eqref{Int3}. These solutions form a $H^{1}$-continuous curve
%\begin{equation}\label{Int9}
%\left\{\varphi_{c,p},~~c\in]0,\delta'_{p}[\right\}\subset H^{4}(\mathbb{R}),
%\end{equation}
of even solitary waves of \eqref{Int1} which are all stable.
\end{Theo}
\begin{Rem}
If in Definition \ref{groundstate}, we replace the requirement $ \beta>0 $ by $ \beta\neq 0 $ then the uniqueness result for the even ground states to \eqref{Int4} with $\mu=1 $ and $ c>0 $ small enough still holds but up to the symmetry $ v\mapsto- v$ in the case $ p $ even.
\end{Rem}

% PROPOSITION DE ANNE
The stability result will follow from a continuity argument together with  a positivity property of  the quadratic form associated with the second Fr\'echet derivative of the action functional at 
 respectively $ \varphi_{1,p,\mu_p} $ and $ \varphi_{1,p,0}$ where $ \mu_p=\frac{2^{2}(p+2)^{2}}
{(p^{2}+4p+8)^{2}}$. Note that to prove this positivity property for $ \varphi_{1,p,\mu_p} $ we  check numerically a sign condition on some $ L^2 $-scalar product linked to $ \varphi_{1,p,\mu_p}$. 

We will apply the following  classical proposition for equations of form the $ \partial_t V'(u)=\partial_x E'(u) $ where $ E $ and $ V $ are conservation laws, $ V$ being quadratic.
\begin{Pro}[see for instance de Bouard \cite{Anne}]\label{Anne}
Let $\phi_{c,p,\mu}\in H^{2}(\mathbb{R})$ be a solution of \eqref{Int3}
traveling with the speed $c>0$, and $\mathcal{L}_{c,p,\mu}$ be
the linearized operator associated to the second derivative of the action functional $ E_{p,\mu}+c V $  at $\phi_{c,p}$, defined by 
$\mathcal{L}_{c,p,\mu}v=\mu\partial^{4}_{x}v-\partial^{2}_{x}v
+cv-\phi_{c,p,\mu}^{p}v$, for all $v\in H^{4}(\mathbb{R})$. 
Assume that the exists $\delta>0$ such that 
\begin{equation}\label{Int10}
\left\langle\mathcal{L}_{c,p,\mu}v,v\right\rangle_{L^{2}}\ge\delta
\left\|v\right\|^{2}_{H^{2}(\mathbb{R})},
\end{equation}
for all $v\in H^{2}(\mathbb{R})$ satisfying the orthogonalities
\begin{equation}\label{Int11}
\left\langle v,\phi_{c,p,\mu}\right\rangle_{L^{2}}=
\left\langle v,\phi'_{c,p,\mu}\right\rangle_{L^{2}}=0.
\end{equation}
Then, $\phi_{c,p,\mu}$ is stable in $H^{2}(\mathbb{R})$. 
\end{Pro}

% ORIGINE DE LA PROPOSITION DE ANNE
The proof of Proposition \ref{Anne} follows 
from the theory developed by Benjamin \cite{MR0338584},
Bona \cite{MR0386438} and Weinstein \cite{MR783974} which relies on the spectral properties of $\mathcal{L}_{c,p,\mu}$
(see Section \ref{Stability-1} for details).

% RESULTATS PARTIELS
For other stability results on solitons  the similar models,
see for instance: Karpman \cite{MR1372681}-\cite{MR1837757},
Dey et al. \cite{MR1427885}, Levandosky \cite{MR1667522}-\cite{MR2332504},
Bridges et al. \cite{MR1946769}, Pava \cite{MR1958041} .

% ORGANISATION DU MANUSCRIPT
This paper is organized as follows. In Section \ref{Stability-1}, 
we prove the assertion  $(i)$ of Theorem \ref{Main-Result} 
and, in Section \ref{Stability-2}, we prove the assertion  $(ii)$
of Theorem \ref{Main-Result}.

% STABILITE AUTOUR DU SOTITON EXPLICITE 
\section{Existence and Stability for the branch crossing the explicit solitary waves  of gKW}
\label{Stability-1}
In this section, we prove the point $(i)$ of Theorem \ref{Main-Result}. First, we observe that $ u $ is a solution to \eqref{Int1} with $ \mu>0 $ if and only if $ u_\mu(t,x)=\mu^{1/p} u(\mu^{3/2} t \sqrt{\mu} x) $ is a solution to \eqref{Int1} 
 with $ \mu=1 $.  In particular, $\varphi $ is a solution to \eqref{Int3}  with $ c=1 $ and $ \mu>0 $ if and only if $ \varphi_{\mu}=\mu^{1/p} \varphi(\sqrt{\mu} \, \cdot) $ is a solution to \eqref{Int3}  with $ c=\mu$ and $ \mu=1 $. This ensures that it is equivalent to prove the assertion $(i) $ of Theorem  \ref{Main-Result} 
 with $ \mu=\mu_p=\frac{2^{2}(p+2)^{2}}{(p^{2}+4p+8)^{2}}$ and $ c $ close to $ 1$. For the numerical checking of the sign of some $ L^2$-scalar product (see Subsection \ref{num}) we prefer to work with this last normalization in this section.

  So, let $p\in\mathbb{N}^{*}$ be fixed and set $ \mu=\mu_p=\frac{2^{2}(p+2)^{2}}{(p^{2}+4p+8)^{2}}$, then  the explicit solitary wave given by \eqref{Int6} travels with the   speed 
   $ c=1$. 
The strategy is as follows. By applying the Implicit Function Theorem, we  will establish for all real $c>0$ sufficiently 
close to $ 1$, the existence of even solitary waves  $\varphi_{c,p,\mu_p}$ of gKW traveling 
with  speed $c$ (see Lemma \ref{TFI}). Next, following Albert \cite{MR1151253}, we use that the hypotheses of Proposition \ref{Anne} are fulfilled for the explicit solitary wave 
 $ \varphi_{1,p,\mu_p} $ whenever a sign condition on some $ L^2 $-scalar product involving  $ \varphi_{1,p,\mu_p} $ is satisfied. This sign condition is checked numerically in Subsection \ref{num}. Finally, arguing by continuity, we prove  that the hypotheses of Proposition \ref{Anne} are still fulfilled for $ c $ close enough to $ 1 $  which leads to the  orbital stability of all the the family 
$\left\{\varphi_{c,p,\mu_p},~c\in]1-\delta_{p}',1+\delta_{p}'[\right\}
\subset H^{4}(\mathbb{R})$, with $0<\delta_{p}'\ll 1$.

% DIVISION DE LA PREUVE EN PLUSIEURS LEMMES
We will divide the proof of Theorem \ref{Main-Result} part $(i)$ 
into several lemmas.  We first study the spectral properties 
of the linearized operator $\mathcal{L}_{1,p,\mu_p}$  associated with the second derivative of the action functional at 
 the explicit soliton $\varphi_{1,p,\mu_p}$
(defined in \eqref{Int6}-\eqref{Int7}). 
 
\subsection{Spectral properties of $\mathcal{L}_{1,p,\mu_p}$}
% PROPRIETES SPECTRALES DE L_c_p
\begin{Lem}[Spectral Properties of $\mathcal{L}_{1,p,\mu_p}$]\label{Spectral}
Let $p\in\mathbb{N}^{*}$. We consider the unbounded operator
$\mathcal{L}_{1,p,\mu_p}:L^{2}(\mathbb{R})\rightarrow L^{2}(\mathbb{R})$, defined by: 
$ u\mapsto \mu_p \partial^{4}_{x}u-\partial^{2}_{x}u+u-\varphi^{p}_{1,p,\mu_p}u$. 
We claim that $\mathcal{L}_{1,p,\mu_p}$ possesses, among others, the following three 
crucial properties:\\
(P1) The essential spectrum of   $\mathcal{L}_{1,p,\mu_p}$ is $ [1,+\infty[ $ ; \\
(P2) $\mathcal{L}_{1,p,\mu_p}$ has  only one negative eigenvalue $\lambda_{1,p,\mu_p}$ which is simple;\\
(P3) The kernel of $\mathcal{L}_{1,p,\mu_p}$ is spanned by $\varphi'_{1,p,\mu_p}$.
\end{Lem}

% PREUVE DES PROPRIETES SPECTRALES 
\textbf{Proof.} 
One can clearly see that $\mathcal{L}_{1,p,\mu_p}$ is linear with domain $H^{4}(\mathbb{R})$, 
self-adjoint  and  closed on $L^{2}(\mathbb{R})$, since $\left\|u\right\|_{H^{4}(\mathbb{R})}
\lesssim \left\|u\right\|_{L^{2}(\mathbb{R})} +\|\mathcal{L}_{1,p,\mu_p} u\|_{L^2} $. Moreover, it is a compact perturbation of $\mu_p \partial^{4}_{x}-\partial^{2}_{x}+1$, 
since $\varphi^{p}_{1,p,\mu_p}$ is smooth and decays exponentially to $0$, and thus, 
its essential spectrum is given by $\sigma_{ess}(\mathcal{L}_{1,p,\mu_p})=[1,+\infty[$.

% THEOREME D'ALBERT
According to Albert \cite{MR1151253} (Theorem $3.2$ and Lemma $10$), 
to get the properties $(P2)$-$(P3)$, it suffices to prove that:
\begin{equation}\label{pf001}
\mathcal{F}(\varphi_{1,p,\mu_p})(\omega)>0,~~\forall\omega\in\mathbb{R},
~~\text{and}~~
\frac{d^{2}}{dw^{2}}\text{log}\mathcal{F}(\varphi^{p}_{1,p,\mu_p})(\omega)<0,
~~\forall\omega\in\mathbb{R}^{*},
\end{equation}
where $\mathcal{F}$ denotes the spatial Fourier transform.
Thanks to the computations done by Magnus and Oberhettinger \cite{MR0232968} (p. $34$),
we first recall that 
\begin{equation}\label{pf02}
\mathcal{F}\left(\text{sech}^{\nu}(\cdot)\right)(\omega)= 2^{\nu-1}(\Gamma(\nu))^{-1}
\left|\Gamma\left(\frac{\nu}{2}+i\frac{\omega}{2}\right)\right|^{2}\ge 0,
~~\forall\omega\in\mathbb{R},~~\nu\in\mathbb{R}_{+},
\end{equation}
where $\Gamma$ denotes the gamma function.
On the other hand, using the usual Fourier transforms table, we establish that
\begin{equation}\label{pf3}
\mathcal{F}\left(\text{sech}^{4}(\cdot/2)\right)(\omega)=
\frac{2^{4}\pi}{3!}\omega(\omega^{2}+1)\text{cosech}(\pi\omega),
~~\forall\omega\in\mathbb{R}
\end{equation}
and thus
\begin{equation}\label{pf5}
\frac{d^{2}}{d\omega^{2}}\log\mathcal{F}\left(\text{sech}^{4}(\cdot/2)\right)(\omega)
=\frac{2^{4}\pi}{3!}\Bigl(  -\frac{1}{\omega^2} +\frac{2(1-\omega^2)}{(1+\omega^2)^2} 
+\frac{\pi^{2}}{\text{sinh}^{2}(\pi\omega)}\Bigr)<0,
~~\forall\omega\in\mathbb{R}^{*}.
\end{equation}
Then, combining \eqref{pf02}-\eqref{pf5}, we get the sufficient condition \eqref{pf001}. 
Thus $\mathcal{L}_{1,p,\mu_p}$ has the spectral properties $(P1)$-$(P3)$.

% INFORMATIONS COMPLEMENTAIRES
Let $\lambda_{1,p,\mu_p}$ be the unique negative (simple) eigenvalue 
of $\mathcal{L}_{1,p,\mu_p}$, and let $\chi_{1,p,\mu_p}$ be the, normalized in $ H^2(\R) $, eigenfunction
 associated with  $\lambda_{1,p,\mu_p}$. One can notice that 
$\left\{\varphi_{1,p,\mu_p},\chi_{1,p,\mu_p}\right\}\subset
(\text{Ker}~\mathcal{L}_{1,p,\mu_p})^{\perp}$
for the usual inner product on $L^{2}(\mathbb{R})$, 
since $\left\langle\chi_{1,p,\mu_p},\varphi'_{1,p,\mu_p}\right\rangle_{L^{2}}
=\lambda^{-1}_{1,p,\mu_p}\left\langle \chi_{1,p,\mu_p},\mathcal{L}_{1,p,\mu_p}
\varphi'_{1,p,\mu_p}\right\rangle_{L^{2}}=0$. Moreover, since $\varphi_{1,p,\mu_p}$
is an even function, then $\mathcal{L}_{1,p,\mu_p}\chi_{1,p,\mu_p}(-\cdot)
=\lambda_{1,p,\mu_p}\chi_{1,p,\mu_p}(-\cdot)$, and by uniqueness 
$\chi_{1,p,\mu_p}(-\cdot)=\chi_{1,p,\mu_p}(\cdot)$. 
This implies that $\chi_{1,p,\mu_p}$ is an even function. 
\hfill $ \square $ \vspace*{2mm}

\subsection{Construction of the $ C^1 $ branch of solitary waves}
We construct now the new solitons in a neighborhood of the explicit
solitary wave $\varphi_{1,p,\mu_p}$. 

% EXISTENCE DE SOLITONS POUR C PROCHE DE c_p
\begin{Lem}[Existence of Solitons $\varphi_{1,p,\mu_p}$ for $c$ close to $1$]\label{TFI}
There exist $ \delta_p>0 $ and $ \tilde{\delta}_p>0 $ such that for any $ c>0 $ with $ |c-1|<\delta_p $, there exists a unique $ H^4(\R) $  even solution 
$ \varphi_{c,p,\mu_p} $ of \eqref{Int3} in the ball of $ H^4(\R) $ centered at $ \varphi_{1,p,\mu_p} $ with radius $ \tilde{\delta}_p>0 $.
Moreover, the function $ c \mapsto  \varphi_{c,p,\mu_p} $  is of class $ C^1 $ from $ ]1-\delta_p, 1+\delta_p[$ into 
$ H^4(\R) $.
\end{Lem}

% PREUVE DE L'EXISTENCE DES SOLITONS 
\textbf{Proof.} We apply of the Implicit Function Theorem,
(see for instance \cite{MR2215276} for a similar application). 
Recall that we  set $H^s_{e}(\R)=\left\{u\in H^{s}(\mathbb{R}):u(-\cdot)=u(\cdot) \right\}$. For $\gamma>0$, we define the 
 map $T:]1-\gamma,1+\gamma[\times H^4_e(\R)\rightarrow L^2_e(\R)$, 
by: $(c,\psi)\mapsto \mu_p \partial^{4}_{x}\psi-\partial^{2}_{x}\psi+c\psi-\frac{1}{p+1}\psi^{p+1}$. 
 $T$ is obviously of class $C^{1}$ and, since  $\varphi_{1,p,\mu_p}$ 
satisfies \eqref{Int3}-\eqref{Int7}, it holds
\begin{equation}\label{exi1}
T(1,\varphi_{1,p,\mu_p})=0\; \text{and}~~
\partial_{\psi}T(1,\varphi_{1,p,\mu_p})=\mathcal{L}_{1,p,\mu_p}\Big\rvert_{H^4_e(\R)},
\end{equation}
where $\mathcal{L}_{1,p,\mu_p}\Bigr\rvert_{H^4_e(\R)}$
denotes the restriction of $\mathcal{L}_{1,p,\mu_p}$ to $H^4_{e}(\R) $.

% ISOMORPHISME 
It is easy to check that the spectrum of  $\mathcal{L}_{1,p,\mu_p}\Bigr\rvert_{H^4_e(\R)}:H^4_e(\R)\rightarrow L^2_e(\R)$
is  contained in the spectrum of $\mathcal{L}_{1,p,\mu_p}$ and   Lemma \ref{Spectral}
(properties (P1) and (P3)) ensures that   $\sigma_{ess}(\mathcal{L}_{1,p,\mu_p})=[1,+\infty[$  
and $\text{Ker}\, \mathcal{L}_{1,p,\mu_p}=\text{span}\,  \{\varphi'_{1,p,\mu_p}\} $. 
 Since   $\varphi'_{1,p,\mu_p}$ is an odd function, we thus infer that  $0\notin\sigma(\mathcal{L}_{1,p,\mu_p}\vert_{H^4_e(\R)})$
  and thus $\mathcal{L}_{1,p,\mu_p}\Bigr\rvert_{H^4_e(\R)}$ is an isomorphism.

% APPLICATION DU THEOREME DES FONCTIONS IMPLICITES
Therefore, according to  the Implicit Function Theorem, there exist $\delta_{p}>0$ and a $C^{1}$ 
map $R:]1-\delta_{p},1+\delta_{p}[\rightarrow H^4_e(\R)$ which is uniquely determined 
such that: $T(c,R(c))=0$ for all $(c,R(c))\in ]1-\delta_{p}, 1+\delta_{p}[
\times B_{H^4_e}(\varphi_{1,p,\mu_p},\delta_{p})$, where $B_{H^4_e}(\varphi_{1,p,\mu_p},\delta_p)$ 
denotes the unit ball of  $ H^4_e(\R)$ centered at  $\varphi_{1,p,\mu_p}$ with radius $ \delta_p$.
 This proves the lemma by setting  $ \varphi_{c,p,\mu_p} =R(c) $ for $ c\in  ]1-\delta_{p}, 1+\delta_{p}[$. \hfill $ \square $ \vspace*{2mm}

% PROPRIETES DES NOUVEAUX SOLITONS DE gKW
\begin{Rem}[Some Important Properties of $\varphi_{c,p}$ for $c\sim 1$]\label{Bootstrap}
\normalfont

Since $\varphi_{c,p,\mu_p}\in H^{4}(\mathbb{R})$, using the equation \eqref{Int3}
of $\varphi_{c,p,\mu_p}$, and a classical bootstrap argument, 
we get that $\varphi_{c,p,\mu_p}\in H^{k}(\mathbb{R})$ for all 
$k\in\mathbb{N}$. The well-known Sobolev embedding of
$H^{k+1}(\mathbb{R})$ into $C^{k}(\mathbb{R})$ leads to 
\begin{equation}\label{lim1}
\left\|\varphi_{c,p,\mu_p}\right\|_{C^{k}([n,n+1])}\le 
C_{S}\left\|\varphi_{c,p,\mu_p}\right\|_{H^{k+1}([n,n+1])}
\rightarrow 0~~\text{as}~~n\rightarrow +\infty,
\end{equation}
and thus $\partial^{k}_{x}\varphi_{c,p,\mu_p}(\pm\infty)=0$
for all $k\in\mathbb{N}$. 

Now, for $c$ close to $1$, in view of the solutions  to the linear asymptotic equation: $\mu_p \partial^{4}_{x}\psi
-\partial^{2}_{x}\psi+c\psi=0$, we infer that  $ \varphi_{c,p,\mu_p} $ and its derivatives satisfy
\begin{equation}\label{lim3}
|\partial^{k}_{x}\varphi_{1,p,\mu_p}(x) |\lesssim  e^{-\sqrt{c}\left|x\right|} \; \text{ for } k\in \{0,1,2,3,4\}.
\end{equation}
\end{Rem}

% FIN DE LA PREUVE DU THEOREME DE STABILITE PARTIE (i)
\subsection{Orbital stability result assuming $\langle\mathcal{L}^{-1}_{1,p,\mu_p}\varphi_{1,p,\mu_p},
\varphi_{1,p,\mu_p}\rangle_{L^{2}}<0$} \label{Stab-gKW}
Recall that, according to Lemma \ref{Spectral}, the operator $\mathcal{L}_{1,p,\mu_p}$ 
possesses the spectral properties $(P1)$-$(P3)$. Then it is well known that if $ w\in H^2(\R) $ satisfies 
$ \langle \mathcal{L}_{1,p,\mu_p} w,w\rangle_{L^2} 	<0 $ then \eqref{Int10} holds for any $ v\in H^2(\R) $ such that $\langle v, \varphi'_{c,p} \rangle_{L^2}
=\langle v,  \mathcal{L}_{1,p,\mu_p} w \rangle_{L^2}=0 $. Now, since $\varphi_{1,p,\mu_p}$ is even, it follows from the last subsection that 
$\varphi_{1,p,\mu_p}\in \text{Im}\, \mathcal{L}_{1,p,\mu_p} $ and applying the above criterium with $ w\in \mathcal{L}^{-1}_{1,p,\mu_p}\varphi_{1,p,\mu_p},$ we obtain  that the hypotheses of Proposition \ref{Anne} are satisfied as soon as
  \begin{equation}\label{main0002}
\mathcal{J}_{p}=\langle\mathcal{L}^{-1}_{1,p,\mu_p}\varphi_{1,p,\mu_p},
\varphi_{1,p,\mu_p}\rangle_{L^{2}}<0 \; .
\end{equation}
Note that the above quantity does not depend on the choice of the element of the preimage of  $\varphi_{1,p,\mu_p} $ since  $ \varphi_{1,p,\mu_p} \in (\ker \,  \mathcal{L}_{1,p,\mu_p})^\perp$. 
In the next subsection we check numerically that this sign condition is fulfilled for $ p\in\{1,2,3,4\} $. In this subsection, we prove the orbital stability result 
 assuming \eqref{main0002}. 

% PRESQUE ORTHOGONALIES AVEC SOLITON EXPLICITE
According to Proposition \ref{Anne} , it suffices to check that for any  $ c$  close enough to $ 1 $ and any  $ v\in H^2(\R ) $ such that 
$\langle v,\varphi_{c,p}\rangle_{L^{2}}=\langle v,\varphi'_{c,p}\rangle_{L^{2}}=0$, there exists $ \delta_{p} >0 $ such that 
\begin{equation}\label{to}
\langle\mathcal{L}_{c,p}v,v \rangle_{L^{2}}\ge \delta_{p} \|v\|_{H^2}^2 \; .
\end{equation}
First, we notice that setting $ \gamma_{c,p} =\|\varphi_{c,p,\mu_p}-\varphi_{1,p,\mu_p} \|_{H^1} $ for $|c-1|<\delta_p $, it holds 
\begin{equation}\label{main3}
\left|\left\langle v,\varphi_{1,p,\mu_p}\right\rangle_{L^{2}}\right|
\le\gamma_{c,p}\left\|v\right\|_{L^{2}(\mathbb{R})}\quad\text{ and }\quad \left|\left\langle v,\varphi'_{1,p,\mu_p}\right\rangle_{L^{2}}\right|
\le\gamma_{c,p}\left\|v\right\|_{L^{2}(\mathbb{R})}.
\end{equation}

% COERCIVITE AVEC VITESSE EXPLICITE
Next, for all $v\in H^{2}(\mathbb{R})$ satisfying the almost orthogonality conditions
\eqref{main3} with $ \gamma_{c,p} $ small enough, let us prove that $\mathcal{L}_{1,p,\mu_p}$ is coercive
in $H^{2}(\mathbb{R})$. We argue as, for instance, C\^ote et al  (\cite {MR3461359}, Lemma $2.6$).  \eqref{main3} clearly leads to 
\begin{equation}\label{posi1}
\left|\left\langle v,\frac{\varphi_{1,p,\mu_p}}{\left\|\varphi_{1,p,\mu_p}
\right\|_{L^{2}(\mathbb{R})}}\right\rangle_{L^{2}}\right|
+
\left|\left\langle v,\frac{\varphi'_{1,p,\mu_p}}
{\left\|\varphi'_{1,p,\mu_p}\right\|_{L^{2}(\mathbb{R})}}\right\rangle_{L^{2}}\right|
\lesssim\gamma_{c,p}\left\|v\right\|_{H^{2}(\mathbb{R})}.
\end{equation}
 Now, we decompose $v$ as follows
\begin{equation}\label{posi2}
v=v_{1}+a_{1}\frac{\varphi_{1,p,\mu_p}}{\left\|\varphi_{1,p,\mu_p}\right\|_{L^{2}}}
+a_{2}\frac{\varphi'_{1,p,\mu_p}}{\left\|\varphi'_{1,p,\mu_p}\right\|_{L^{2}}}
=v_{1}+v_{2},
\end{equation}
with $\left\langle v_{1},\varphi_{1,p,\mu_p}\right\rangle_{L^{2}}=
\left\langle v_{1},\varphi'_{1,p,\mu_p}\right\rangle_{L^{2}}=0$. 
Then, combining \eqref{posi1} and \eqref{posi2}, we infer that
\begin{equation}\label{posi3}
\left|a_{1}\right|+\left|a_{2}\right|\lesssim\gamma_{c,p}
\left\|v\right\|_{H^{2}(\mathbb{R})}.
\end{equation}
Moreover, for $k=0,1,2$, using that 
$\left\langle\partial^{k}_{x}\varphi_{1,p,\mu_p},\partial^{k}_{x}\varphi'_{1,p,\mu_p}
\right\rangle_{L^{2}}=0$, it clearly holds
\begin{equation}\label{posi4}
\left\|\partial^{k}_{x}v_{2}\right\|^{2}_{L^{2}(\mathbb{R})}=
a^{2}_{1}\frac{\left\|\partial^{k}_{x}\varphi_{1,p,\mu_p}
\right\|^{2}_{L^{2}(\mathbb{R})}}
{\left\|\varphi_{1,p,\mu_p}\right\|^{2}_{L^{2}(\mathbb{R})}}
+
a^{2}_{2}\frac{\left\|\partial^{k}_{x}\varphi'_{1,p,\mu_p}
\right\|^{2}_{L^{2}(\mathbb{R})}}
{\left\|\varphi'_{1,p,\mu_p}\right\|^{2}_{L^{2}(\mathbb{R})}}.
\end{equation}
Then, combining \eqref{posi2}-\eqref{posi4},
we deduce that there exists $ 0<\varepsilon_p<1 $ such that  
\begin{equation}\label{posi5}
\left\|v_{1}\right\|_{H^{2}(\mathbb{R})}\sim 
\left\|v\right\|_{H^{2}(\mathbb{R})} \quad \text{as soon as} \; \gamma_{c,p}<\varepsilon_p\; .
\end{equation}
Next, we compute 
\begin{equation}\label{posi6}
\left\langle\mathcal{L}_{1,p,\mu_p}v,v\right\rangle_{L^{2}}
=\left\langle\mathcal{L}_{1,p,\mu_p}v_{1},v_{1}\right\rangle_{L^{2}}
+\left\langle\mathcal{L}_{1,p,\mu_p}v_{2},v_{2}\right\rangle_{L^{2}}
+2\left\langle\mathcal{L}_{1,p,\mu_p}v_{1},v_{2}\right\rangle_{L^{2}}.
\end{equation}
Note that, by construction, $v_{1}$ is orthogonal in $ L^2(\R) $ with $ \varphi_{1,p,\mu_p} $ and $ \varphi'_{1,p,\mu_p}$. Therefore,  according to the discussion in the beginning of this Subsection,  there exists $\delta_{0}>0$ such that 
\begin{equation}\label{posi7}
\left\langle\mathcal{L}_{1,p,\mu_p}v_{1},v_{1}\right\rangle_{L^{2}}\ge 
\delta_{0}\left\|v_{1}\right\|^{2}_{H^{2}(\mathbb{R})}\gtrsim \delta_{0}
\left\|v\right\|^{2}_{H^{2}(\mathbb{R})}.
\end{equation}
On the other hand, using \eqref{posi3}-\eqref{posi5}, it holds
\begin{equation}\label{posi8}
\left|\left\langle\mathcal{L}_{1,p,\mu_p}v_{2},v_{2}\right\rangle_{L^{2}}\right|
\lesssim\left|a_{1}\right|^{2}+\left|a_{2}\right|^{2}
\lesssim (\gamma_{c,p})^2\left\|v\right\|^{2}_{H^{2}(\mathbb{R})}
\end{equation}
and 
\begin{equation}\label{posi9}
\left|\left\langle\mathcal{L}_{1,p,\mu_p}v_{1},v_{2}\right\rangle_{L^{2}}\right|
\lesssim\left\|v_{1}\right\|_{H^{2}(\mathbb{R})}
\left\|v_{2}\right\|_{H^{2}(\mathbb{R})}
\lesssim\left\|v\right\|_{H^{2}(\mathbb{R})}
(\left|a_{1}\right|+\left|a_{2}\right|)
\lesssim\gamma_{c,p}\left\|v\right\|^{2}_{H^{2}(\mathbb{R})}.
\end{equation}
Therefore, combining \eqref{posi6}-\eqref{posi9}, we deduce that there exists $ \delta'_p >0$ and $ \varepsilon_p>0 $ such that if 
 $ \|\varphi_{c,p}-\varphi_{c_p,p}\|_{L^2}< \varepsilon_p $ then 
\begin{equation}\label{posi10}
\left\langle\mathcal{L}_{1,p,\mu_p}v,v\right\rangle_{L^{2}}
\ge \delta'_{p}\left\|v\right\|^{2}_{H^{2}(\mathbb{R})}\, .
\end{equation}
 \eqref{to} follows immediately by noticing that 

\begin{align}\label{posi11}
\left\langle\mathcal{L}_{c,p,\mu_p}v,v\right\rangle_{L^{2}}&
=\left\langle\mathcal{L}_{1,p,\mu_p}v,v\right\rangle_{L^{2}}
+(c-1)\left\|v\right\|^{2}_{L^{2}(\mathbb{R})}
-\left\langle v^{2},\varphi^{p}_{c,p,\mu_p}-\varphi^{p}_{1,p,\mu_p}
\right\rangle_{L^{2}}
\nonumber\\
&\ge \delta'_{p}\left\|v\right\|^{2}_{H^{2}(\mathbb{R})}
-(c-1+K \varepsilon_p)\left\|v\right\|^{2}_{L^{2}(\mathbb{R})}
\nonumber\\
&\ge \frac{\delta'_{p}}{2}\left\|v\right\|^{2}_{H^{2}(\mathbb{R})},
\end{align}
as soon as $|c-1+K\varepsilon_p|< \delta'_p / 2$.  But the continuity of the branch of the $ \varphi_{c,p,\mu_p} $  in $ H^4(\R )$ ensures that this is true as soon as 
 $ c$ is close enough to $ 1$ and we are done.

% ETAPE 1
\subsection{Numerical checking of the sign condition on   $\langle\mathcal{L}^{-1}_{1,p,\mu_p}\varphi_{1,p,\mu_p},
\varphi_{1,p,\mu_p}\rangle_{L^{2}}$}\label{num}
\subsubsection{Numerical computing of 
$\mathcal{J}_{1}=\left\langle\mathcal{L}^{-1}_{1,1,\mu_1}\varphi_{1,1,\mu_1},
\varphi_{1,1,\mu_1}\right\rangle_{L^{2}}$.}

% CONFIGURATION DES PARAMETRES ET REMARQUES
Let $\rho_{1} \in\mathcal{L}^{-1}_{1,1,\mu_1}\varphi_{1,1}$
Note that, since $\mathcal{L}_{1,1,\mu_1}$ is a self-adjoint
operator on $L^{2}(\mathbb{R})$, the value of $\mathcal{J}_{1}$ 
does not depend on the choice of the element $\rho_{1}$. Moreover, using that 
$\mathcal{L}_{1,1,\mu_1}\varphi_{1,1,\mu_1}=-\frac{1}{2}\varphi^{2}_{1,1,\mu_1}$, we observe that
\begin{equation}\label{pf01}
\left\langle\rho_{1},\varphi^{2}_{1,1,\mu_1}\right\rangle_{L^{2}}
=-2\|\varphi_{1,1,\mu_1}\|^{2}_{L^{2}(\mathbb{R})}<0.
\end{equation}
From \eqref{pf01} and the exponential decay properties of $\varphi_{1,1,\mu_1}$
and $\rho_{1}$ (see \eqref{lim3}), we infer that $\rho_{1}$ takes necessarily
some negative values near the origin. Then, we can expect that $\mathcal{J}_{1}$ is negative.

\begin{figure}[ht]
\centering
\subfloat[$\varphi_{1,1,\mu_1}(x)$ and $\rho_{1}(x)$ profiles.]
{\includegraphics[width=8cm, height=4cm]{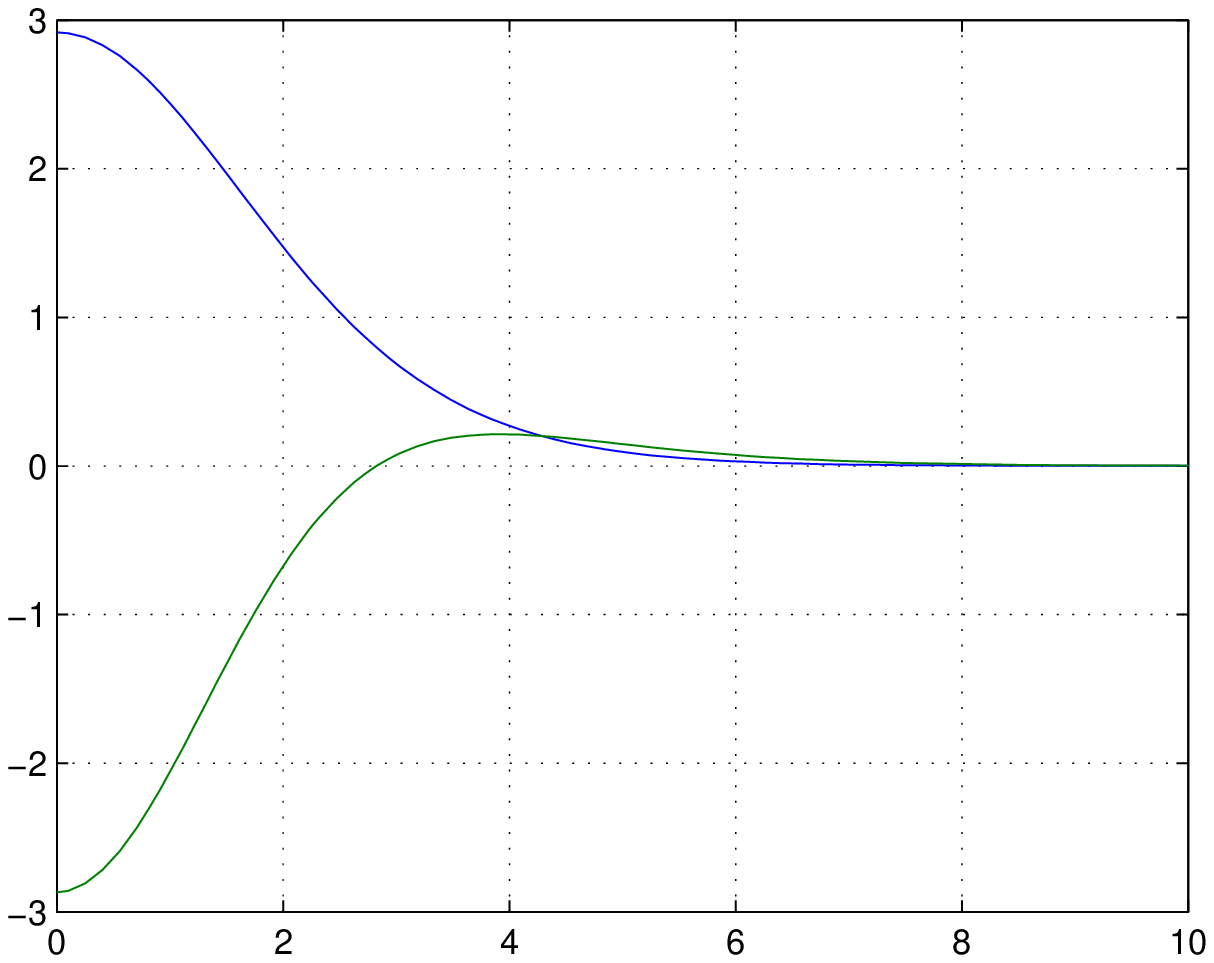}\label{fig4}}
\subfloat[$\int_{0}^{x}\rho_{1}(y)\varphi_{1,1,\mu_1}(y)$ profile.]
{\includegraphics[width=8cm, height=4cm]{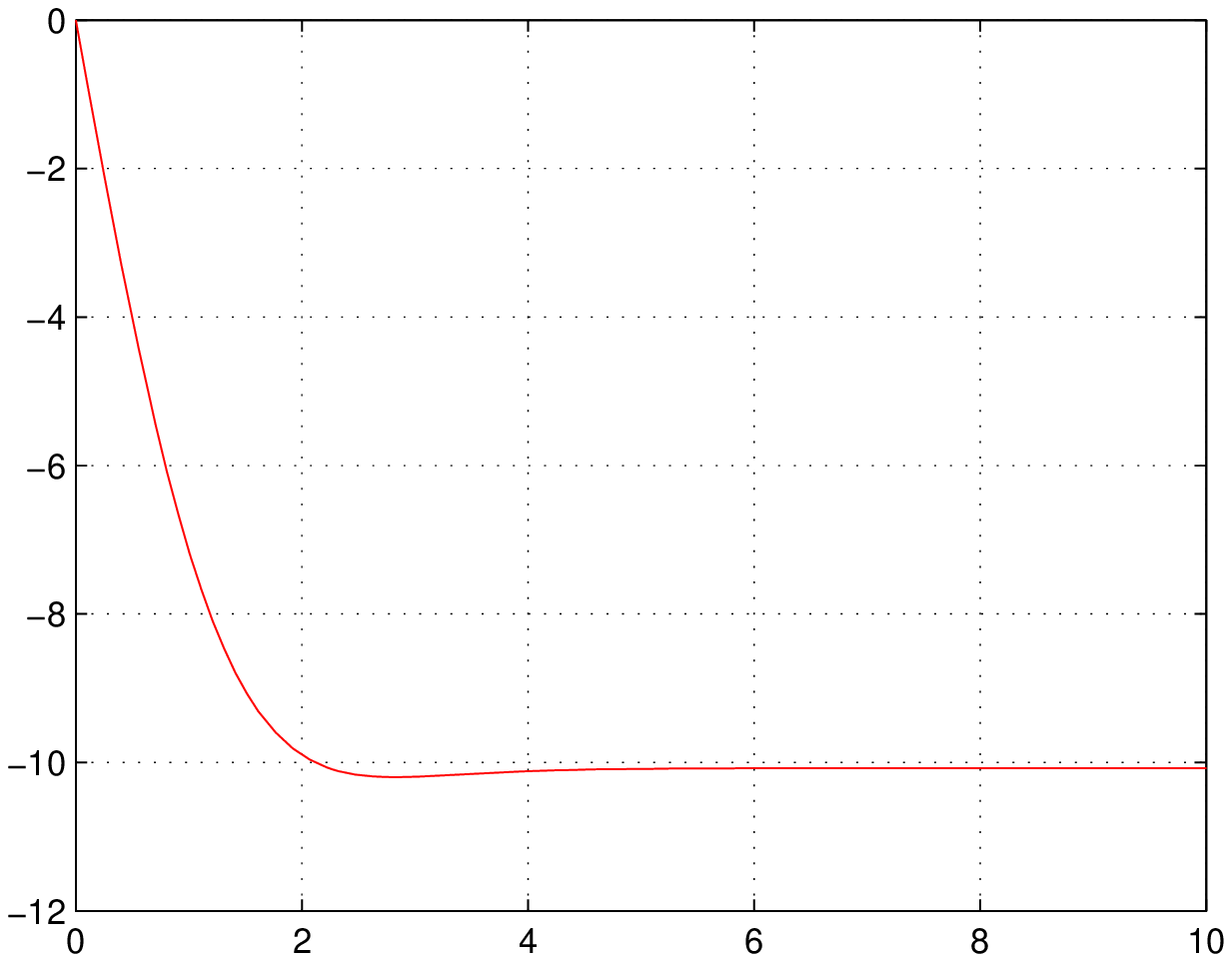}\label{fig5}}
 \\
\subfloat[$\partial^{k}_{x}\varphi_{1,1}(x),~k=0,\ldots,4$, profiles.]
{\includegraphics[width=8cm, height=4cm]{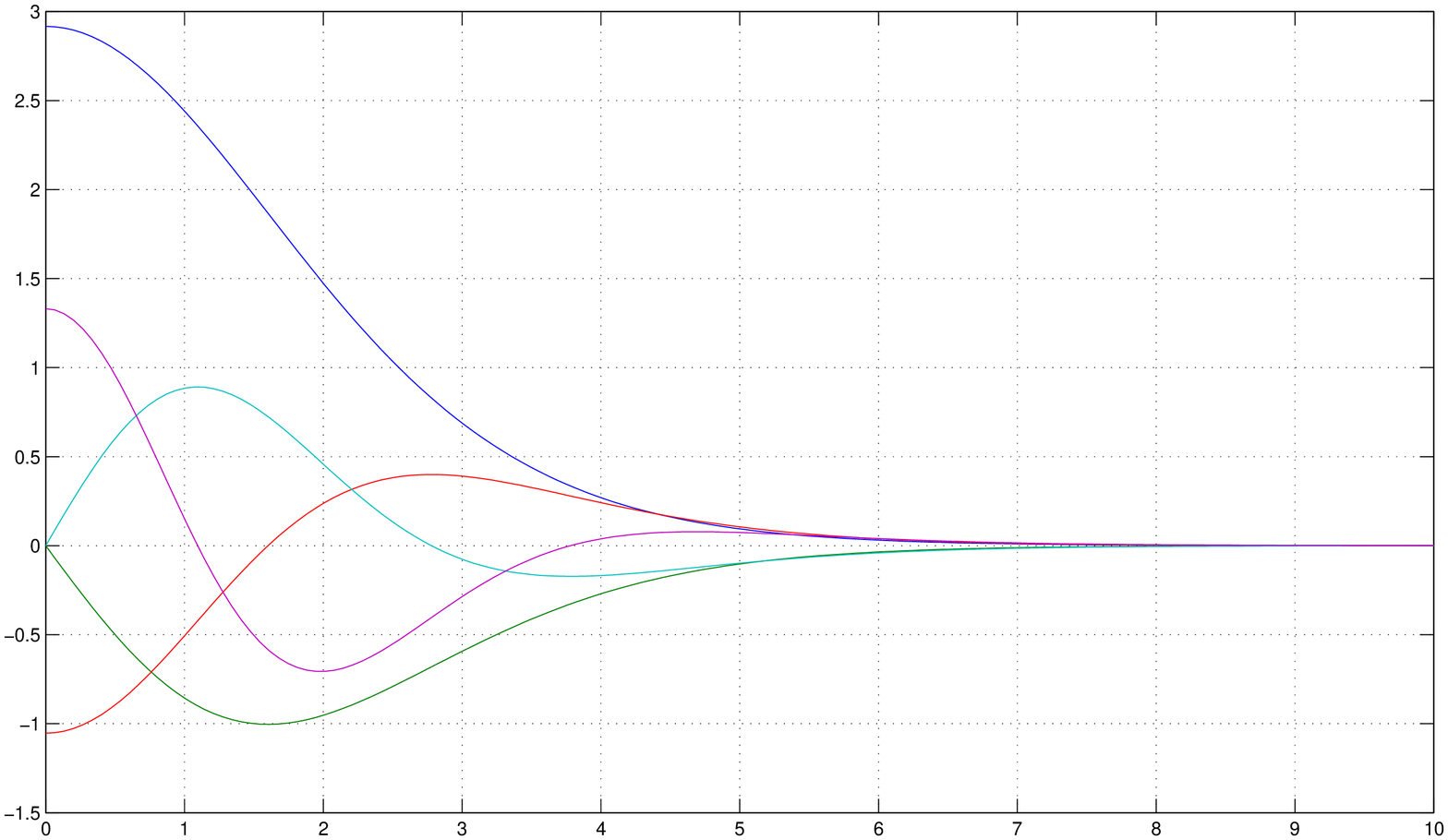}\label{fig6}}
\subfloat[$\partial^{k}_{x}\rho_{1}(x),~k=0,\ldots,4$, profiles.]
{\includegraphics[width=8cm, height=4cm]{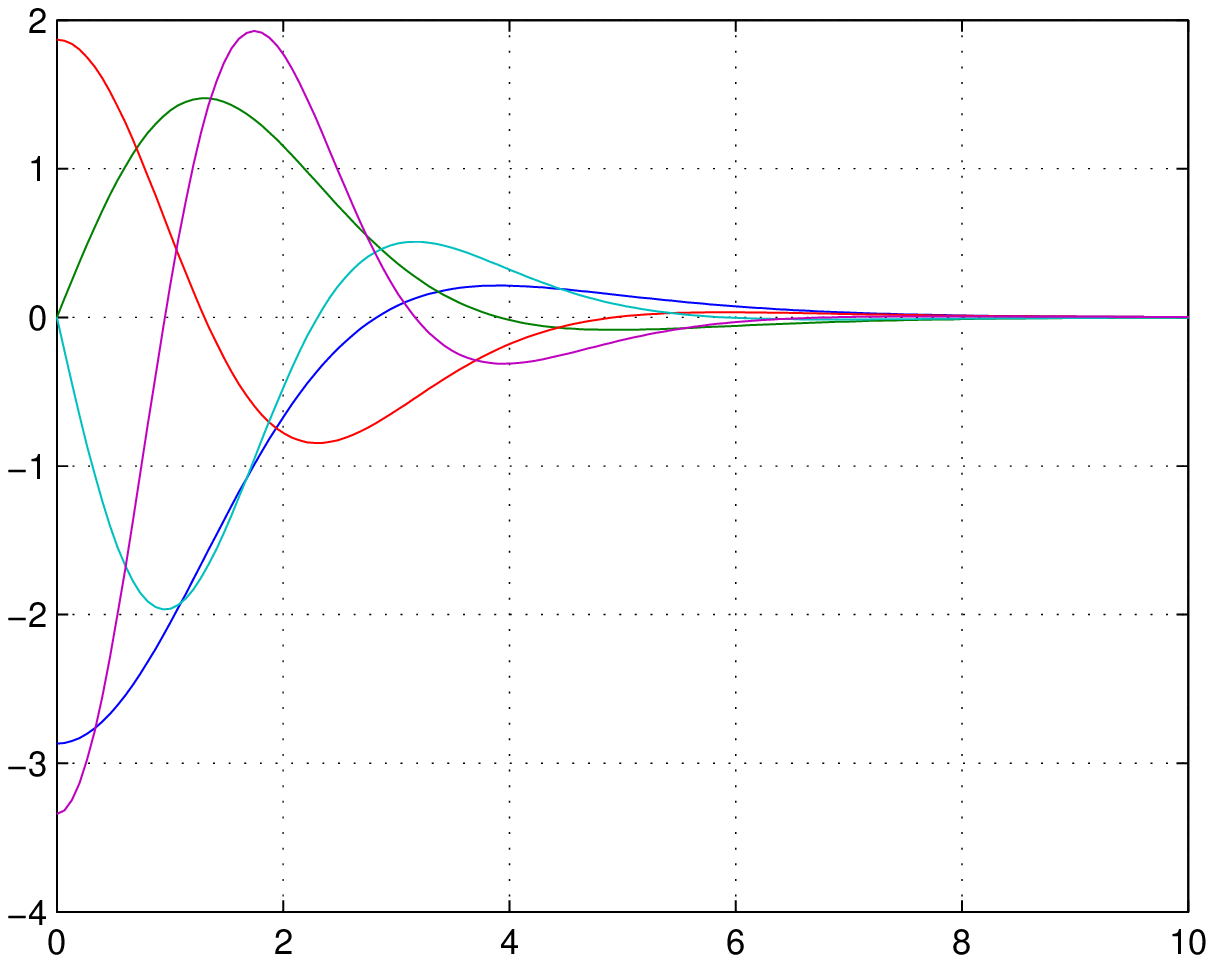}\label{fig7}}
\caption{Numerical computation of  $\mathcal{J}_{1}$}
\end{figure}

% METHODE NUMERIQUE 
The goal is to solve numerically the following equation:
\begin{equation}\label{pf6}
\partial^{4}_{x}\rho_{1}=\frac{13^{2}}{6^{2}}
\left(\partial^{2}_{x}\rho_{1}-\rho_{1}+\varphi_{1,1,\mu_1}\rho_{1}+\varphi_{1,1,\mu_1}
\right).
\end{equation}
First, we remark that $\rho_{1}(\cdot)$ and $\rho_{1}(-\cdot)$ are both solutions 
of \eqref{pf6}, since $\varphi_{1,1,\mu_1}$ is an even function. Then, 
from \eqref{pf01} we deduce that $\rho_{1}$ is also an even function.
Moreover, we get that $\partial^{2k}_{x}\rho_{1}$ 
is even and $\partial^{2k+1}_{x}\rho_{1}$ is odd for all $k\in\mathbb{N}$.
We can restrict ourself to study equation \eqref{pf6} on $\mathbb{R}_{+}$.
We fix the domain to $[0,r_{\max}]$, with $0<r_{\max}<+\infty$,
and we rewrite equation \eqref{pf6} as a system of four first-order ODE. 
We set: $\rho_{1,1}=y_{1}$, $y'_{1}=y_{2}$, $y'_{2}=y_{3}$, 
$y'_{3}=y_{4}$ and $y'_{4}=\frac{13^{2}}{6^{2}}(y_{3}-y_{1}
+\varphi_{1,1,\mu_1}y_{1}+\varphi_{1,1,\mu_1})$. Next, we choose the Robin type boundary
conditions at point $r_{\max}$: $y_{1}(r_{\max})
+y_{2}(r_{\max})=0$ and $y_{3}(r_{\max})+y_{4}(r_{\max})=0$,
and at point $0$ we use the symmetry of $\rho_{1,1}$: $y_{2}(0)=0$ and $y_{4}(0)=0$.
Finally, in our numerical scheme, we take into account 
the exponential decay properties: $\partial^{k}_{x}\varphi_{1,1,\mu_1}(x)\sim 
\partial^{k}_{x}\rho_{1}(x)\sim e^{-x}$, for $x\sim r_{\max}$,
and for all $k\in\mathbb{N}$. Therefore, computing with the MATLAB solver, the result is: 
$\left\langle\rho_{1,1},\varphi_{1,1}\right\rangle_{L^{2}([0,r_{\max}])}\approx-10.0787<0$ 
(see Fig. \ref{fig4}-\ref{fig5}).

% POSITIVITE DE L'OPERATEUR KAWAHARA

% STABILITE DANS LE CAS gKDV CRITIQUE
\subsubsection{Numerical values of
\texorpdfstring{$\mathcal{J}_{p}=\left\langle\mathcal{L}^{-1}_{1,p,\mu_p}\varphi_{1,p,\mu_p},
\varphi_{1,p,\mu_p}\right\rangle_{L^{2}}$}{TEXT}
for \texorpdfstring{$p\in\left\{2,3,4,5\right\}$}{TEXT}} \label{Critic}
\begin{figure}[ht]
\centering
\subfloat[$\varphi_{1,4,\mu_4}(x)$ and $\rho_4\in\mathcal{L}_{1,4,\mu_4}^{-1}\varphi_{1,4}(x)$ 
 profiles.]
{\includegraphics[width=8cm, height=4cm]{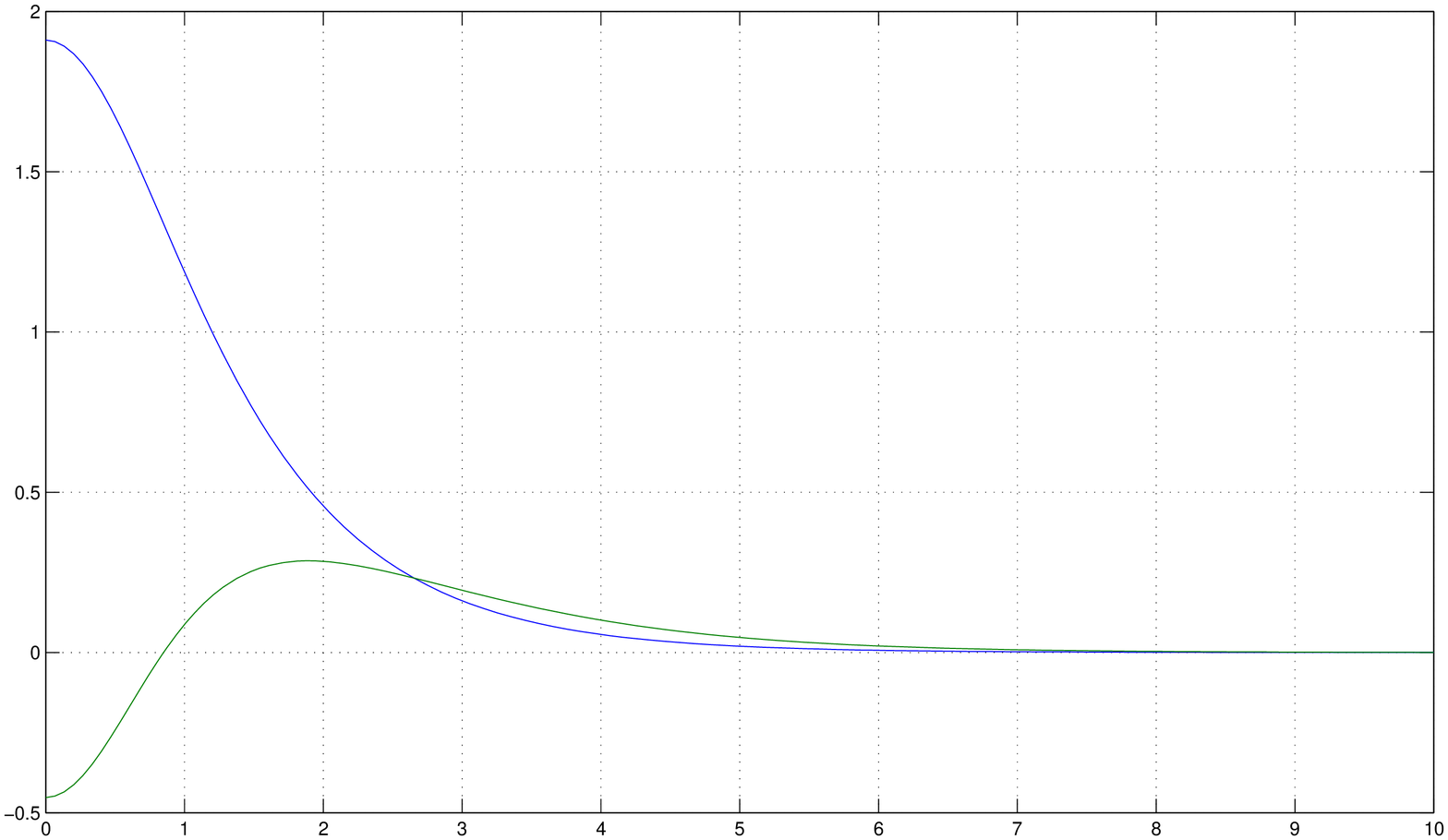}\label{fig12}} 
\subfloat[$\int_{0}^{x}(\mathcal{L}_{1,4,\mu_4}^{-1}\varphi_{1,4,\mu_4}(y))\varphi_{1,4,\mu_4}(y)$ 
 profile.]
{\includegraphics[width=8cm, height=4cm]{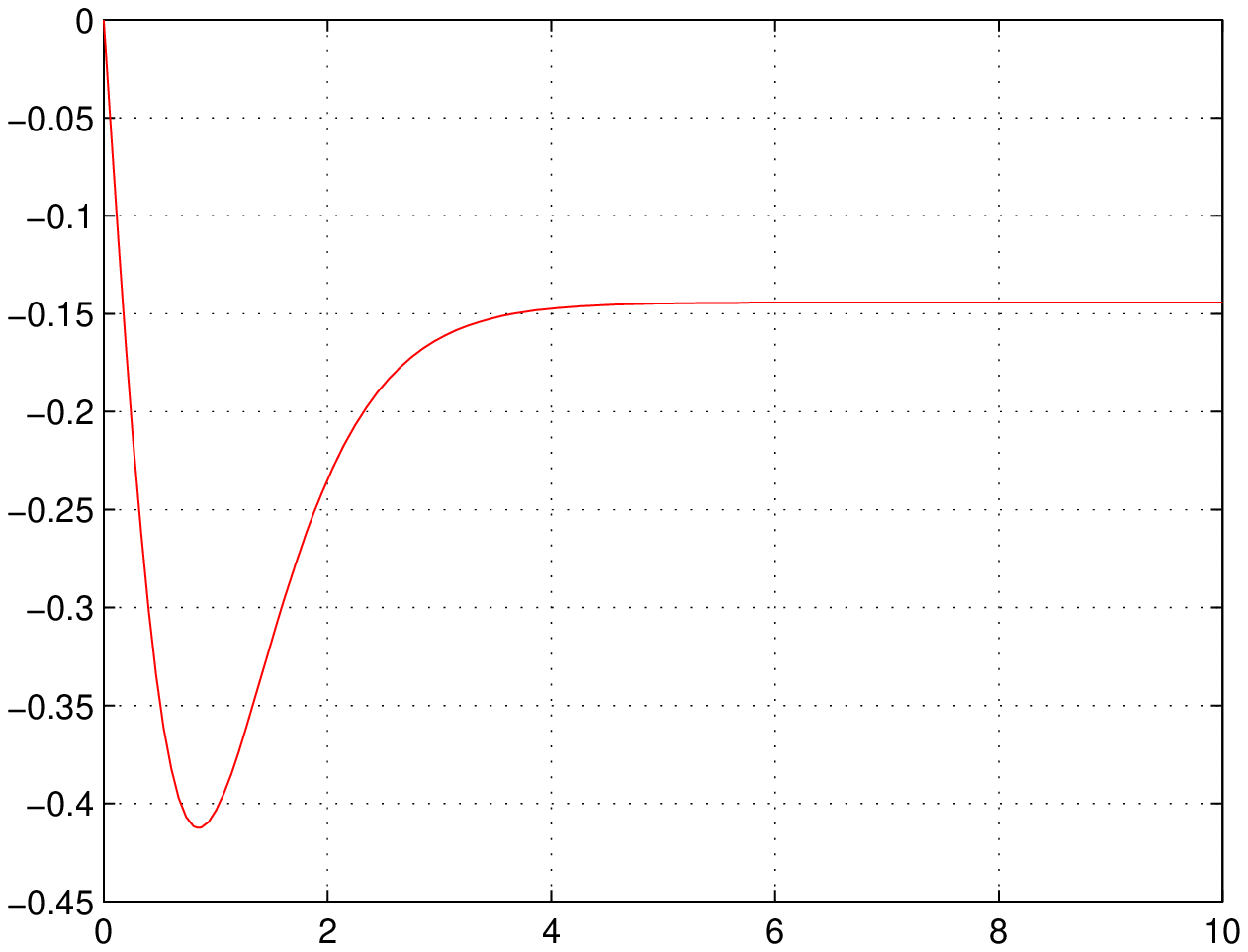}\label{fig13}}
\\
\subfloat[$\varphi_{1,5,\mu_5}(x)$ and $\rho_5=\mathcal{L}_{1,5,\mu_5}^{-1}\varphi_{1,5,\mu_5}(x)$  profiles.]
{\includegraphics[width=8cm, height=4cm]{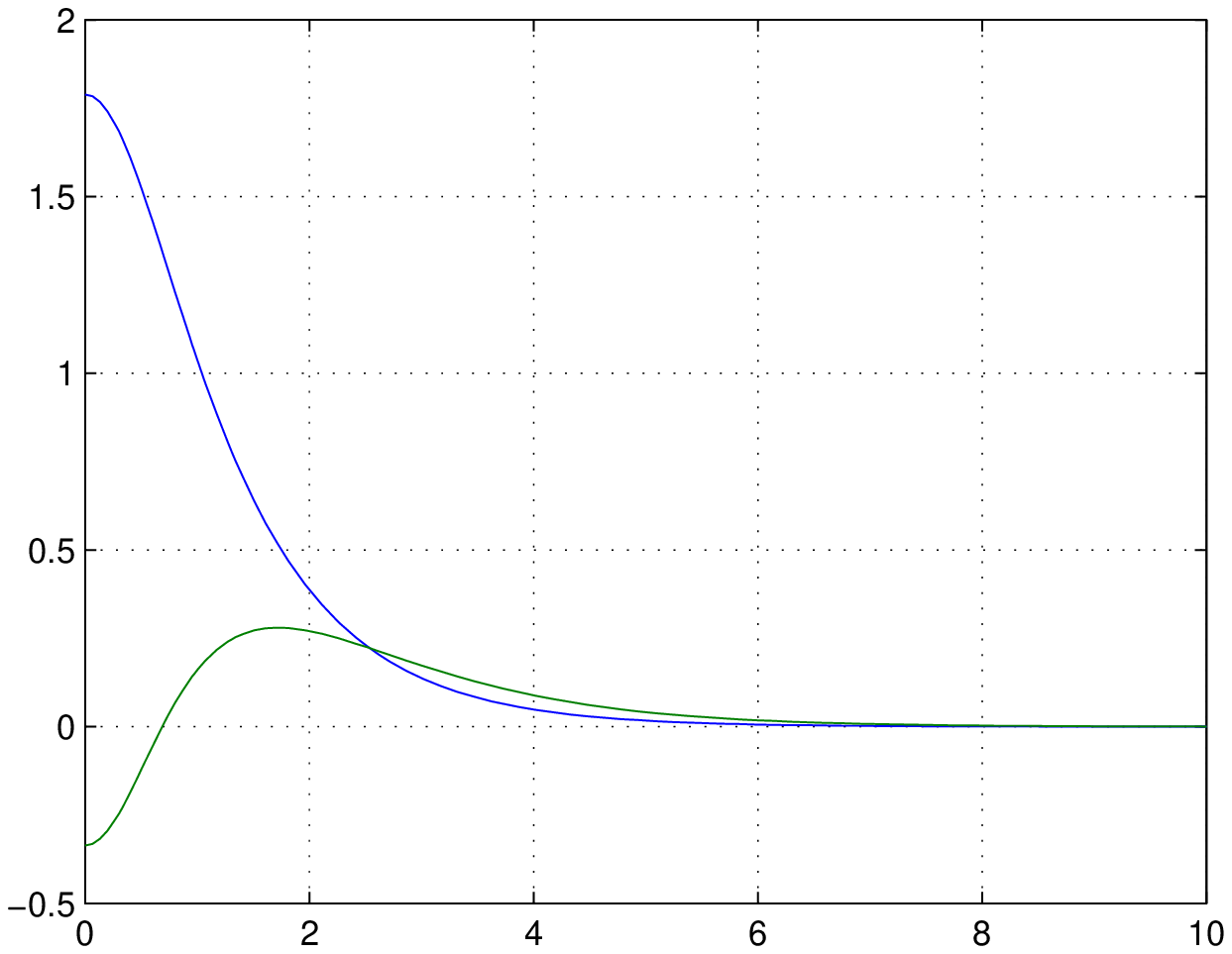}\label{fig14}} 
\subfloat[$\int_{0}^{x}(\mathcal{L}_{1,5,\mu_5}^{-1}\varphi_{1,5,\mu_5}(y))\varphi_{1,5,\mu_5}(y)$  profile.]
{\includegraphics[width=8cm, height=4cm]{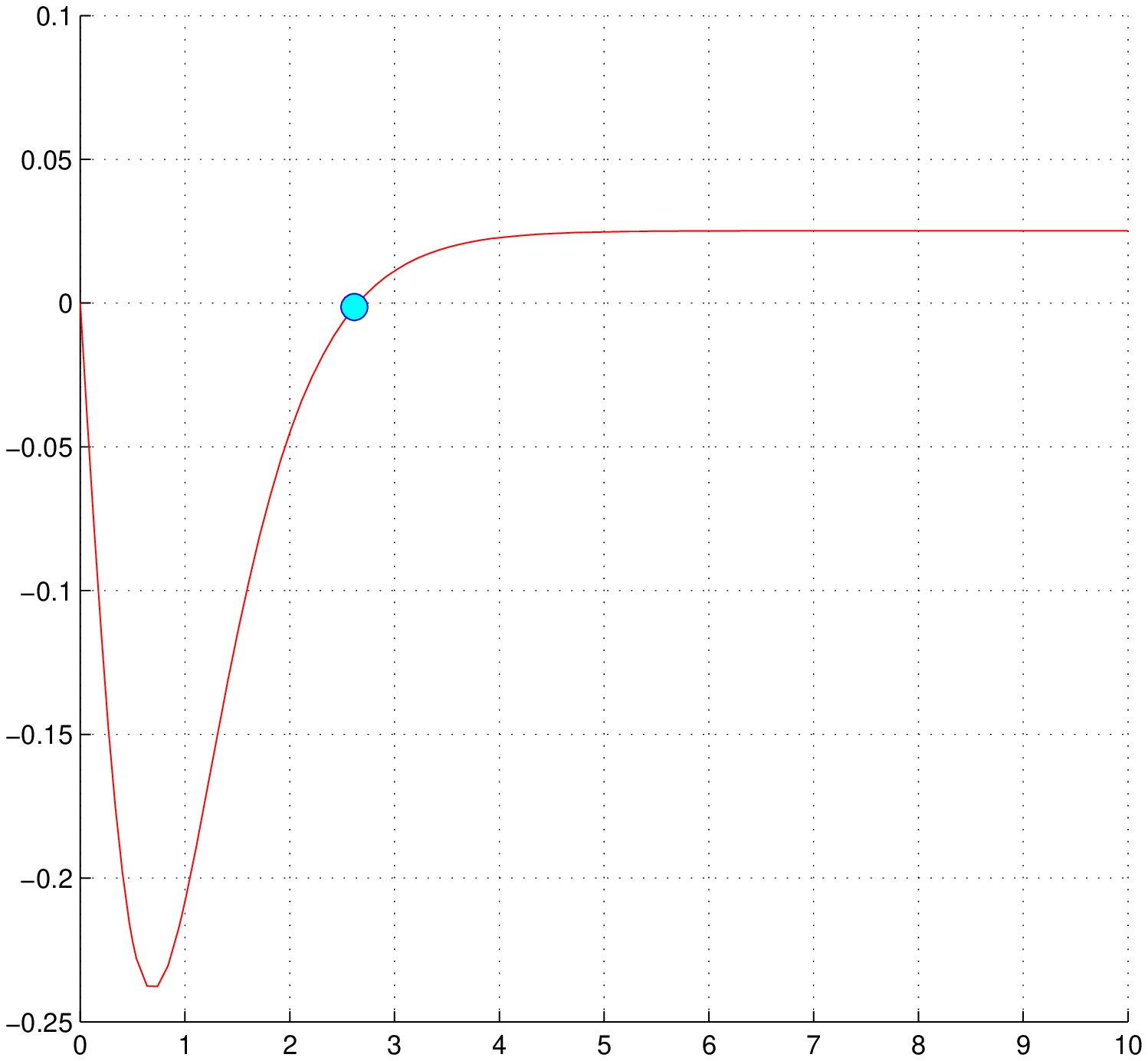}\label{fig15}}
\caption{Numerical computation of $\mathcal{J}_{p}$  for $ p=4,5$.}
\end{figure}

Proceeding as in Subsection \ref{Stab-gKW}, Step $1$, for $p\in\{2,3,4,5\}$, 
we compute respectively that:
$\mathcal{J}_{2}\approx-1.9325<0$, $\mathcal{J}_{3}\approx-0.5649<0$,
$\mathcal{J}_{4}\approx-0.1443<0$ (see Fig. \ref{fig12}-\ref{fig13}) and 
$\mathcal{J}_{5}\approx 0.0252>0$ (see Fig. \ref{fig14}-\ref{fig15}). 
More precisely, we lose the sufficient condition of stability exactly 
for $p_{crit}=4.84$. We observe that the same critical
value appears in the stability analysis done by Bridges 
et al. \cite{MR1946769} (see Section $8$, p. $209$-$210$) 
on the similar model, using the Evans functions approach.
Finally, we point out that thanks to the term $\partial^{5}_{x}u$,
the family $\left\{\varphi_{c,p,1},~c\in]c_{p}-\delta_{p},c_{p}
+\delta_{p}[\right\}\subset H^{4}(\mathbb{R})$, 
with $0<\delta_{p}\ll 1$, of even solitons of gKW
remains stable in the critical gKdV case $p=4$.

% CONDITION SUFFISANTE DE STABILITE POUR P=4 ET p=5

% gKW VERS gKdV
\section{Existence and stability of slow solitons of gKW}
\label{Stability-2} 
In this section, we prove the point $(ii)$ of Theorem \ref{Main-Result}. 
The idea is as follows. Around the functional equation
$E'_{\mu}(\varphi_{1,p,\mu})+V'(\varphi_{c,p,\mu})=0$ 
(equivalent to equation \eqref{Int3}), we will construct a minimization
problem on the line. By applying the Concentration-Compactness
Principe \cite{MR778970},
we will prove the existence of even solitary waves  
$\varphi_{1,p,\mu}$ of gKW, traveling with the speed $1$.
Next, we will establish the strong convergence in $H^{1}(\mathbb{R})$
of the family $\left\{\varphi_{1,p,\mu},~0<\mu\ll 1\right\}$ 
to the solitons $\varphi_{1,p}=\varphi_{1,p,0}$
of gKdV as $\mu$ tends to $0^{+}$
(see Subsection \ref{Mini}).  
% RAPPELS SOLITONS DE gKdV 
Once again, we recall that the solitons $\varphi_{c,p,0}$ 
of gKdV satisfy equation \eqref{Int3}
with $\mu=0$, i.e.:
\begin{equation}\label{preli0}
-\partial^{2}_{x}\varphi_{c,p,0}(x)+c\varphi_{c,p,0}(x)
=\frac{1}{p+1}\varphi^{p+1}_{c,p,0}(x),
~~\forall x\in\mathbb{R}.
\end{equation}
They are unique, up to the translations, and the transformations: 
$\varphi_{c,p}\mapsto-\varphi_{c,p}$ if $p$ is even. Moreover, 
they are explicitly given by:
\begin{equation}\label{Min01}
\varphi_{c,p}(x)=\left[\frac{(p+1)(p+2)c}{2}\right]^{1/p}\text{sech}^{2/p}
\left[\frac{p\sqrt{c}}{2}x\right],~~\forall x\in\mathbb{R}.
\end{equation}

% STABILITE POUR LES SOLITONS DE PETITE VITESSE
Finally, using the well-known spectral properties 
of the linearized operateur 
$\mathcal{L}_{c,p,0}
=-\partial^{2}_{x}(\cdot)+c(\cdot)-\varphi^{p}_{c,p,0}$
of gKdV around the solitons $\varphi_{c,p,0}$, and arguing
as in Subsection \ref{Stab-gKW}, we will prove the stability
in $H^{2}(\mathbb{R})$ of the family $\left\{\varphi_{1,p,\mu},
~0<\mu\ll 1\right\}$ for the sub-critical gKdV nonlinearity
$p\in\left\{1,2,3\right\}$. Note that we also obtain a uniqueness result on the solution to our minimization problem.  After an  appropriate rescaling,
this enables us to construct a continuous branch of  even solitary waves  of gKW traveling 
with low speeds that are all  orbitally stable .

% PROBLEME DE MINIMISATION 
\subsection{Existence and limit of  even ground states of gKW as 
\texorpdfstring{$\mu$}{TEXT} tends to 
\texorpdfstring{$0^{+}$}{TEXT}}\label{Mini}

In this subsection, we first prove the existence of even solitary waves 
$\varphi_{1,p,\mu}$ of gKW for $\mu>0$, by solving a minimization
problem on the line. In this step we argue as in  Levandosky \cite{MR1667522} (Theorem $2.3$).
Second, we study the strong convergence
in $H^{1}(\mathbb{R})$ of the family $\left\{\varphi_{1,p,\mu},
~0<\mu\ll 1\right\}$ to the explicit solitons $\varphi_{1,p}$
of gKdV (given in \eqref{Min01}) as $\mu$ goes to $0^{+}$
(see Lemma \ref{Exist-Conv}). Moreover, we establish
the uniqueness of $\varphi_{1,p,\mu}$ 
for $0<\mu\ll1$  and $ p\in \{1,2,3\} $(see Lemma \ref{Uniq}).

% RESULTAT D'EXISTENCE ET  DE CONVERGENCE
\begin{Lem}[Existence and limit of $\varphi_{1,p,\mu}$]\label{Exist-Conv}
Let $p\in\mathbb{N}^{*}$ be fixed. There exists a family 
$\left\{\varphi_{1,p,\mu},~0<\mu\ll 1\right\}\subset H^{4}(\mathbb{R})$ 
of even solutions to  \eqref{Int3} with $ c=1$, such that 
\begin{equation}\label{Min04}
\lim_{\mu\to 0^{+}}\left\|\varphi_{1,p,\mu}
-\varphi_{1,p,0}\right\|_{H^{1}(\mathbb{R})}=0.
\end{equation}
\end{Lem}

% PREUVE
\textbf{Proof.} We split the proof in five steps.

% ETAPE 1: PROBLEME DE MINIMIZATION SUR R
\textit{Step 1.} Minimization problem on $\mathbb{R}$.

We argue similarly as Levandosky \cite{MR1667522} (Theorem $2.3$).
First, we set  $H^2_e(\R)=\left\{\psi\in H^{2}(\mathbb{R}):~\psi(-\cdot)
=\psi(\cdot)\right\}$. Clearly, $H^2_e(\R)$, endowed with the scalar product of  $H^{2}(\mathbb{R})$, is a Hilbert space since we have a continuous embedding of
$H^{2}(\mathbb{R})$ into $C^{1}(\mathbb{R})$. For all $\psi\in H^2_e(\R)$,
we define the functional associated with the linear part
of  \eqref{Int3}  by:
\begin{equation}\label{min1}
I_{\mu}(\psi)=I_{1,\mu}(\psi)=\int_{\mathbb{R}}\left[\frac{\mu}{2}(\partial^{2}_{x}\psi)^{2}
+\frac{1}{2}(\partial_{x}\psi)^{2}+\frac{1}{2}\psi^{2}\right],
\end{equation}
and the functional associated with the nonlinear part of 
\eqref{Int3} by:
\begin{equation}\label{min2}
K_{p}(\psi)=\frac{1}{(p+1)(p+2)}\int_{\mathbb{R}}\psi^{p+2}.
\end{equation}
One can see that $I_{\mu}(\cdot)$ is coercive in $H^{2}(\mathbb{R})$,
since for all $\psi\in H^{2}(\mathbb{R})$, we have
\begin{equation}\label{min1-0}
I_{\mu}(\psi)\ge\frac{\min\{\mu,1\}}{2}\left\|\psi\right\|^{2}_{H^{2}(\mathbb{R})}.
\end{equation}
Also, one can check that $K_{p}(\cdot)$ is locally Lipschitz in $L^{p+2}(\mathbb{R})$. 
Indeed, for all $(\psi,\phi)\in L^{p+2}(\mathbb{R})\times L^{p+2}(\mathbb{R})$, 
applying the H\"older inequality, we have
\begin{align}\label{min1-1}
\left|K_{p}(\psi)-K_{p}(\phi)\right|&=\frac{1}{(p+1)(p+2)}
\left|\int_{\mathbb{R}}(\psi^{p+2}-\phi^{p+2})\right|
\nonumber\\
&\le\frac{1}{(p+1)(p+2)}\left\|\psi-\phi\right\|_{L^{p+2}(\mathbb{R})}
\sum_{k=0}^{p+1}\left\||\psi|^{p-k+1}|\phi|^{k}\right\|_{L^{\frac{p+2}{p+1}}(\mathbb{R})},
\end{align}
and applying the Young inequality, we have 
\begin{equation}\label{min1-2}
\left\||\psi|^{p-k+1}|\phi|^{k}\right\|^{\frac{p+2}{p+1}}_{L^{\frac{p+2}{p+1}}(\mathbb{R})}
\le\frac{p-k+1}{p+1}\left\|\psi\right\|^{p+2}_{L^{p+2}(\mathbb{R})}
+\frac{k}{p+1}\left\|\phi\right\|^{p+2}_{L^{p+2}(\mathbb{R})}.
\end{equation}
Then combining \eqref{min1-1} and \eqref{min1-2}, we obtain
\begin{equation}\label{min1-3}
\left|K_{p}(\psi)-K_{p}(\phi)\right|\le 
C(\left\|\psi\right\|_{L^{p+2}(\mathbb{R})},\left\|\phi\right\|_{L^{p+2}(\mathbb{R})})
\left\|\psi-\phi\right\|_{L^{p+2}(\mathbb{R})}.
\end{equation}

% PROBLEME DE MINIMIZATION 
Now, we will study the minimization of $I_{\mu}(\cdot)$ in $H^2_e(\R)$ 
subject to the constraint $K_{p}(\cdot)=K_{p}(\varphi_{1,p,0})=\beta_{p}$ where $ \varphi_{1,p,0} $ is defined in \eqref{Int0003}.
We chose this constraint to be sure that the family of minimizer
(depending on the parameter $\mu$) will converge to the solitons $\varphi_{1,p,0}$
of gKdV (in a sense that we specify in Step $5$). For $ \beta>0 $, we define:
\begin{equation}\label{min4}
S^{\beta}_{p,\mu}=S^{\beta}_{1,p,\mu}=\inf\left\{I_{\mu}(\psi):~\psi\in H^2_e(\R),
~K_{p}(\psi)=\beta\right\}.
\end{equation}
Since $\varphi_{1,p}\in H^2_e(\R)$, and $I_{\mu}(\cdot)$ is coercive in $H^2_e(\R)$,
one can see that $S_{p,\mu}^{\beta_p}>0$. On the other hand, using the definition 
of $I_{\mu}(\cdot)$ and $K_{p}(\cdot)$, one can check that
\begin{equation}\label{min4-1}
S^{\beta}_{p,\mu}=\left(\frac{\beta}{\beta_{p}}\right)^{\frac{2}{p+2}}
S^{\beta_{p}}_{p,\mu},~~\forall\beta\in\mathbb{R}^{*}_{+},
\end{equation}
and then we infer the crucial sub-additivity property:
\begin{equation}\label{min4-2}
S^{\beta}_{p,\mu}+S^{\beta_{p}-\beta}_{p,\mu}>S^{\beta_p}_{p,\mu},
~~\forall\beta\in(0,\beta_{p}).
\end{equation}

% RESOLUTION DU PROBLEME DE MINIMISATION SUR R
Let us solve the constraint minimization problem \eqref{min4} with $ \beta=\beta_p $.
We take $(\psi_{k})_{k\ge1}\subset H$ a minimizing 
sequence of the problem \eqref{min4}, i.e., for all $k\in\mathbb{N}^{*}$, 
we have
\begin{equation}\label{min5}
K_{p}(\psi_{k})=\beta_{p}~~\text{and}~~
\lim_{k\to+\infty}\left|I_{\mu}(\psi_{k})-S_{p,\mu}\right|=0.
\end{equation}
We claim that there exist a sub-sequence $(\psi_{j})_{j\ge 1}$ and
a function $\psi_{\mu}\in H^2_e(\R)$, such that $\psi_{j}\rightarrow\psi_{\mu}$
strongly in $H^{2}(\mathbb{R})$ as $j\rightarrow+\infty$. Moreover,
$\psi_{\mu}$ is a minimizer of the problem \eqref{min4}.

% EXTRACTION D'UNE SOUS-SUITE
From the convergence of the energy \eqref{min5}, we deduce that 
there exists $\varepsilon_{0}>0$, such that 
\begin{equation}\label{min5-1}
0<I_{c,\mu}(\psi_{k})<S_{c,p,\mu}+\varepsilon,
\end{equation}
with $0<\varepsilon<\varepsilon_{0}\ll 1$, and this implies that the sequence 
$(\psi_{k})_{k\ge 1}$ is bounded in $H^{2}(\mathbb{R})$ (since $I_{c,\mu}(\cdot)$ 
is equivalent to $H^{2}(\mathbb{R})$). Since $H$ is reflexive,
there exist a sub-sequence $(\psi_{j})_{j\ge 1}\subset H$ 
and a function $\psi_{\mu}\in H$ such that the following hold:
\begin{equation}\label{min5-2}
\psi_{j}\rightharpoonup\psi_{\mu}~~\text{weakly}~~\text{in}~~H^{2}(\mathbb{R}),
\end{equation}
\begin{equation}\label{min5-2-1}
\psi_{j}\rightarrow\psi_{\mu}~~\text{strongly}~~\text{in}~~H^{1}_{loc}(\mathbb{R}),
\end{equation}
\begin{equation}\label{min5-2-0}
\psi_{j}\rightarrow\psi_{\mu}~~\text{a.e. on}~~\mathbb{R},
\end{equation}
and 
\begin{equation}\label{min5-2-2}
\left\|\psi_{\mu}\right\|^{2}_{H^{2}(\mathbb{R})}
\le\liminf_{j\to+\infty}\left\|\psi_{j}\right\|^{2}_{H^{2}(\mathbb{R})}.
\end{equation}
Now,  we define a sequence of positive and even functions:
\begin{equation}\label{min5-3}
\Phi_{k}=\left|\partial^{2}_{x}\psi_{k}\right|^{2}
+\left|\psi_{k}\right|^{2},~~\forall k\in\mathbb{N}^{*}.
\end{equation}
From \eqref{min5-1}, one can see that $(\Phi_{k})_{k\ge 1}$ 
is bounded in $L^{1}(\mathbb{R})$.
After extracting a sub-sequence, we may assume that 
$\lim_{k\to +\infty}\int_{\mathbb{R}}\Phi_{k}=L<+\infty$. By normalizing,
we may assume further that $\int_{\mathbb{R}}\Phi_{k}=L$ for all 
$k\in\mathbb{N}^{*}$. Then, by the Concentration-Compactness
Lemma \cite{MR1667522}, there are three possibilities:

% CONCENTRATION
$(a)$ Compactness: there exists $(y_{k})_{k\ge 1}\subset\mathbb{R}$,
such that for all $\varepsilon>0$, there exists $R_{\varepsilon}>0$, 
such that for all $k\in\mathbb{N}^{*}$,
\begin{equation}\label{min5-4}
\int_{\left|x-y_{k}\right|\le R_{\varepsilon}}\Phi_{k}\ge L-\varepsilon
\Leftrightarrow\int_{\left|x-y_{k}\right|\ge R_{\varepsilon}}
\Phi_{k}\le\varepsilon.
\end{equation}

% EVANESSANCE 
$(b)$ Vanishing: for every $R>0$,
\begin{equation}\label{min5-5}
\lim_{k\to+\infty}\sup_{y\in\mathbb{R}}\int_{\left|x-y\right|\le R}
\Phi_{k}=0.
\end{equation}

% DICHOTOMIE
$(c)$ Dichotomy: there exists $l\in(0,L)$, such that for all $\varepsilon>0$, 
there exist $R>0$, $R_{k}\rightarrow+\infty$, $(y_{k})_{k\ge 1}$ and $k_{0}$, 
such that for all $k\ge k_{0}$, 
\begin{equation}\label{min6-1}
\left|\int_{\left|x-y_{k}\right|\le R}\Phi_{k}-l\right|\le\varepsilon
~~\text{and}~~
\left|\int_{R<\left|x-y_{k}\right|<R_{k}}\Phi_{k}\right|\le\varepsilon. 
\end{equation}

% EXISTENCE DU MINIMISEUR 
First, let us assume that $(a)$ holds, and we will prove that $\psi_{\mu}$
is a minimizer of the problem \eqref{min4}. One can remark
that $|y_{k}|\le R_{L/2}$ for all $k\in\mathbb{N}^{*}$, otherwise,
there exists $k_{0}$ such that $\{|x-y_{k_{0}}|\le R_{L/2}\}\subset\mathbb{R}^{+}$.
Then, using that $\Phi_{k_{0}}$ is even, we get a contradiction:
\begin{equation}\label{min6-1-0}
\int_{\mathbb{R}}\Phi_{k_{0}}>2\int_{|x-y_{k_{0}}|<R_{L/2}}\Phi_{k_{0}}>L.
\end{equation}
Thus, we can assume that $y_{k}=0$ for all $k\in\mathbb{N}^{*}$ 
in hypothesis \eqref{min5-4} by taking as radius $R_{L/2}+R_{\varepsilon}$.
Now, by the compact embedding of $H^{1}_{loc}(\mathbb{R})$
into $L^{p+2}_{loc}(\mathbb{R})$, and using 
\eqref{min5-2-1}, we obtain
\begin{equation}\label{min6-2}
\psi_{j}\rightarrow\psi_{\mu}~~\text{strongly in}~~L^{p+2}_{loc}(\mathbb{R}).
\end{equation}
On the other hand, using that $\left\|\psi_{j}\right\|_{L^{p+2}(\mathbb{R})}
\le \frac{1}{2^{p/2(p+2)}}\left\|\psi_{j}\right\|_{H^{1}(\mathbb{R})}$,
and the compactness of $(\Phi_{j})_{j\ge 1}$, one can easily
check that the sequence $(\left|\psi_{j}\right|^{p+2})_{j\ge 1}$
 is also compact, i.e. 
 $$
 \lim_{R\to +\infty} \sup_{j\in\N} \int_{|x-y_j|>R} |\psi_j|^{p+2} =0\; .
 $$
This clearly leads to 
\begin{equation}\label{min6-3}
\psi_{j}\rightarrow\psi_{\mu}~~\text{strongly in}~~L^{p+2}(\mathbb{R}).
\end{equation}
 Finally, by lower semi-continuity 
\eqref{min5-2-2}, we get $I_{\mu}(\psi_{\mu})\le S^{\beta_{p}}_{p,\mu}$,
and from the fact that $K_{p}(\cdot)$
is locally Lipschitz in $L^{p+2}(\mathbb{R})$ 
(see \eqref{min1-3}), and using \eqref{min6-3},
we get $K_{p}(\psi_{\mu})=\beta_{p}$. 
Therefore, $\psi_{\mu}$ is a minimizer 
of the problem \eqref{min4}.

% EXCLUSION DE L'EVANESSANCE
Next, we suppose that $(b)$ holds and prove that this leads
to a contradiction. One can easy estimate that for all $y\in\mathbb{R}$,
\begin{equation}\label{min6-6}
\Bigl| \int_{\left|x-y\right|\le 1}\psi^{p+2}_{k} \Bigr| \lesssim
\left(\int_{\left|x-y\right|\le 1}\Phi_{k}\right)^{\frac{p+2}{2}}.
\end{equation}
Then, using $(b)$, for any fixed $ \varepsilon>0 $, there exists $ k_\varepsilon \ge 0 $ such that for $ k\ge k_\varepsilon$,
\begin{equation}\label{min7}
\Bigl| \int_{\left|x-y\right|\le 1}\psi^{p+2}_{k} \Bigr| \lesssim
\varepsilon^{\frac{p+1}{2}}.
\end{equation}
Now, multiplying \eqref{min7} by $\frac{1}{(p+1)(p+2)}$ 
and summing over intervals
centered at even integers, we obtain
\begin{equation}\label{min7-1}
K_{p}(\psi_{k})\lesssim\varepsilon^{\frac{p-1}{2}} \to 0 \text{ as } \varepsilon \to 0 \; .
\end{equation}
This contradicts the fact that $K_{p}(\psi_{k})=\beta_{p}$
for all $k\in\mathbb{N}^{*}$.

% EXCLUSION DE LA DICHOTOMIE
At last, we assume that $(c)$ holds, and we will rule out
this possibility. We define two smooth non negative test-functions
($C^{\infty}(\mathbb{R})$)  as follows:
\begin{equation}\label{min7-2}
  \xi_{1}(x)=\left\{
    \begin{aligned}
     &1,~~\left|x\right|\le 1,\\
     &0,~~\left|x\right|>2,\\
    \end{aligned}
  \right.
~~\text{and}~~
 \xi_{2}(x)=\left\{
    \begin{aligned}
     &1,~~\left|x\right|\ge 1,\\
     &0,~~\left|x\right|\le\frac{1}{2}.\\
    \end{aligned}
  \right.
\end{equation}
Next, we construct two sequences $\psi_{k,1}(x)
=\xi_{1}(\left|x-y_{k}\right|/R)\psi_{k}(x)$
and $\psi_{k,2}(x)=\xi_{2}(\left|x-y_{k}\right|/R_{k})\psi_{k}(x)$, 
for all $x\in\mathbb{R}$ and $k\in\mathbb{N}^{*}$. 
Using the assumption $(c)$ and the continuous embedding from $ H^1(\R) $ into $ L^{p+2}(\R) $, one can check that for $k\ge k_{0}$,
\begin{equation}\label{min7-3}
I_{\mu}(\psi_{k})=I_{\mu}(\psi_{k,1})
+I_{\mu}(\psi_{k,2})+O(\varepsilon)
\end{equation}
and 
\begin{equation}\label{min7-4}
K_{p}(\psi_{k})=K_{p}(\psi_{k,1})+K_{p}(\psi_{k,2})+O(\varepsilon).
\end{equation}
Since $(\psi_{k})_{k\ge 1}$ is uniformly bounded in $H^{2}(\mathbb{R})$
and $\left\|\psi_{k,i}\right\|_{H^{2}(\mathbb{R})}
\le\left\|\psi_{k}\right\|_{H^{2}(\mathbb{R})}$, for $i=1,2$, 
it follows that $(\psi_{k,i})_{k\ge 1}$, for $i=1,2$, 
are also bounded uniformly
in $H^{2}(\mathbb{R})$. Hence $K_{p}(\psi_{k,i})$, for $i=1,2$,
are bounded and we may pass to a sub-sequence to define 
$\beta_{i}(\varepsilon)=\lim_{k\to+\infty}K_{p}(\psi_{k,i})$, for $i=1,2$. 
Since $\beta_{i}(\varepsilon)$, for $i=1,2$, are bounded uniformly, 
we can choose a sequence $\varepsilon_{j}\rightarrow 0$
such that $\beta_{p,i}=\lim_{j\to+\infty}\beta_{i}(\varepsilon_{j})<\infty$,
for $i=1,2$. We clearly get $\beta_{p,1}+\beta_{p,2}=\beta_{p}$, 
and there are three cases to consider.

% CAS 1
Case $1$: if $\beta_{p,1}\in(0,\beta_{p})$. Using \eqref{min4-1}
and \eqref{min7-3}, we compute that
\begin{align}\label{min7-5}
I_{\mu}(\psi_{k})&\ge S^{K_{p}(\psi_{k,1})}_{p,\mu}+S^{K_{p}(\psi_{k,2})}_{p,\mu}
+O(\varepsilon_{j})
\nonumber\\
&=\left[\left(\frac{K_{p}(\psi_{k,1})}{\beta_{p}}\right)^{\frac{2}{p+2}}
+\left(\frac{K_{p}(\psi_{k,2})}{\beta_{p}}\right)^{\frac{2}{p+2}}\right]
S^{\beta_{p}}_{p,\mu}+O(\varepsilon_{j}).
\end{align}
Letting $k\rightarrow+\infty$, and using that $(\psi_{k})_{k\ge 1}$ is a minimizing 
sequence, we obtain
\begin{equation}\label{min7-6}
S^{\beta_{p}}_{p,\mu}\ge \left[\left(\frac{\beta_{p,1}}{\beta_{p}}\right)^{\frac{2}{p+2}}
+\left(\frac{\beta_{p,2}}{\beta_{p}}\right)^{\frac{2}{p+2}}\right]
S^{\beta_{p}}_{p,\mu}+O(\varepsilon_{j}).
\end{equation}
Letting $j\rightarrow+\infty$, we get a contradiction:
\begin{equation}\label{min7-7}
S^{\beta_{p}}_{p,\mu}\ge \left[\left(\frac{\beta_{p,1}}{\beta_{p}}\right)^{\frac{2}{p+2}}
+\left(\frac{\beta_{p,2}}{\beta_{p}}\right)^{\frac{2}{p+2}}\right]
S^{\beta_{p}}_{p,\mu}
>S^{\beta_{p}}_{p,\mu}.
\end{equation}

% CAS 2
Case $2$: if $\beta_{1}=0$ (or equivalently $\beta_{1}=\beta_{c,p}$). 
Using the coercivity of $I_{\mu}(\cdot)$ in $H^{2}(\mathbb{R})$,
and the assumption $(c)$, we have
\begin{align}\label{min8}
I_{\mu}(\psi_{k,1})&=\int_{\mathbb{R}}
\left(\left|\partial^{2}_{x}\psi_{k,1}\right|^{2}
+\left|\psi_{k,1}\right|^{2}\right)
\nonumber\\
&=\int_{\left|x-y_{k}\right|\le 2R}
\left(\left|\partial^{2}_{x}\psi_{k}\right|^{2}
+\left|\psi_{k}\right|^{2}\right)+O(\varepsilon_{j})
\nonumber\\
&=l+O(\varepsilon_{j}).
\end{align}
Thus, combining \eqref{min7-3} and \eqref{min8}, we get
\begin{equation}\label{min8-1}
I_{\mu}(\psi_{k})\ge l+
\left(\frac{K_{p}(\psi_{k,2})}{\beta_{p}}\right)^{\frac{2}{p+2}}
S^{\beta_{p}}_{p,\mu}+O(\varepsilon_{j}).
\end{equation}
Letting $k,j\rightarrow+\infty$, we get a contradiction:
\begin{equation}\label{min8-2}
S^{\beta_{p}}_{p,\mu}\ge l+
S^{\beta_{p}}_{p,\mu}>S^{\beta_{p}}_{p,\mu}.
\end{equation}

% CAS 3
Case $3$: if $\beta_{1}>\beta_{c,p}$ (or equivalently $\beta_{1}<0$).
Using the non negativity of $I_{\mu}(\cdot)$, we have 
\begin{equation}\label{min8-3}
I_{\mu}(\psi_{k})\ge I_{\mu}(\psi_{k,1})+O(\varepsilon_{j})
\ge\left(\frac{K_{p}(\psi_{k,1})}{\beta_{p}}\right)^{\frac{2}{p+2}}
S^{\beta_{p}}_{p,\mu}+O(\varepsilon_{j}).
\end{equation}
As previously, letting $k,j\rightarrow+\infty$, 
once again, we get a contradiction:
\begin{equation}\label{min8-4}
S^{\beta_{p}}_{p,\mu}\ge \left(\frac{\beta_{1}}
{\beta_{p}}\right)^{\frac{2}{p+2}}
S^{\beta_{p}}_{p,\mu}>S^{\beta_{p}}_{p,\mu}.
\end{equation}
Therefore, we deduce that the sequence $(\psi_{k})_{k\ge 1}$ is compact.

% EQUATION D'EULER-LAGRANDE
\textit{Step 2.} Euler-Lagrange equation.

At this stage, we can write the Euler-Lagrange equation
related to the minimization problem \eqref{min4}.  Since $ \psi\mapsto I_p(\psi) $ and $\psi\mapsto K_p(\psi) $ are obviously of class $ C^1 $ in $H^2(\R) $, the minimizer $\psi_{p,\mu}$
satisfies the following functional equation:
\begin{equation}\label{min10}
\mu\langle\partial^{2}_{x}\psi_{p,\mu},\partial^{2}_{x}\phi\rangle_{L^{2}}+
\langle\partial_{x}\psi_{p,\mu},\partial_{x}\phi\rangle_{L^{2}}+
\langle\psi_{p,\mu},\phi\rangle_{L^{2}}=\frac{\alpha_{p,\mu}}{p+1}
\langle\psi^{p+1}_{p,\mu},\phi\rangle_{L^{2}},
~~\forall\phi\in H^2_e(\R),
\end{equation}
where $\alpha_{p,\mu}\in\mathbb{R}^{*}$ is a Lagrange multiplier. 
It is worth noticing that,   equation \eqref{min10} also holds with any 
odd test-function $\phi\in H^{2}(\mathbb{R})$, since then both members do cancel. 
Therefore $\psi_{p,\mu}$ satisfies equation \eqref{min10} 
for all $\phi\in H^{2}(\mathbb{R})$. 
Now, by a standard bootstrap argument, we get that 
$\psi_{p,\mu}\in H^{s}(\mathbb{R})$ and $\partial^{s}_{x}
\psi_{p,\mu}(\pm\infty)=0$
for all $s\in\mathbb{N}$. Therefore $\psi_{p,\mu}$ is a 
strong solution of \eqref{min10}, i.e.:
\begin{equation}\label{min010}
\mu\partial^{4}_{x}\psi_{p,\mu}(x)-\partial^{2}_{x}
\psi_{p,\mu}(x)+\psi_{p,\mu}(x)
=\frac{\alpha_{p,\mu}}{p+1}\psi^{p+1}_{\mu}(x),
~~\forall x\in\mathbb{R}.
\end{equation}

% DECROISSANCE EXPONENTIELLE 
For $0<\mu<\frac{1}{4}$, let us set:
\begin{equation}\label{min0001}
s_{1}=\pm\frac{\sqrt{1-\sqrt{1-4\mu }}}{\sqrt{2\mu}}
 ~~\text{and}~~
s_{2}=\pm\frac{\sqrt{1+\sqrt{1-4\mu }}}{\sqrt{2\mu}}.
\end{equation}
One can check that $e^{s_{1}x}$ and $e^{s_{2}x}$ are solutions of 
the linear asymptotic equation:
$\mu\partial^{4}_{x}\psi-\partial^{2}_{x}\psi
+\psi=0$.
Then, we infer the following exponential decay properties:
\begin{equation}\label{min0002}
| \partial^{s}_{x}\psi_{p,\mu}(x)|\lesssim  e^{-\left|x\right|},
 ~~\text{for}~~
\left|x\right|\gg 1,
~~s\in\mathbb{N},
~~\text{and}~~
0<\mu\ll 1.
\end{equation}

% CONTROLE DU MULTIPLICATEUR DE LAGRANGE
\textit{Step 3.} Refined estimate on Lagrange multipliers.

We will prove that the family $\left\{\alpha_{p,\mu},~0<\mu< 1\right\}$ 
of Lagrange multipliers is bounded. Multiplying \eqref{min010} by $\psi_{p,\mu}$, 
and integrating on space, we obtain
\begin{equation}\label{Min12}
\alpha_{p,\mu}=\frac{2I_{\mu}(\psi_{p,\mu})}{(p+2)\beta_{p}}=\frac{2S^{\beta_p}_{p,\mu}}{(p+2)\beta_{p}} \le \frac{2I_{\mu}(\varphi_{1,p,0})}{(p+2)\beta_{p}}
\end{equation}
where $ \varphi_{1,p,0} $ is the soliton of the gKdV equation given by \eqref{Int0003}. From the uniqueness of $ \varphi_{1,p,0} $ (up to the symmetries of the equation) 
 and the definition of $ \beta_{p} $, we know that $ \varphi_{1,p} $ is the unique solution to the constraint minimization problem $ S^{\beta_p}_{p,0}  $ 
  which forces 
  \begin{equation}\label{Min112}
 \frac{2I_{0}(\varphi_{1,p,0})}{(p+2)\beta_{p}}=1 \; 
 \end{equation}
 and thus
 \begin{equation}\label{Min13}
\alpha_{p,\mu}\le 1+O(\mu) \; .
\end{equation}
 This last estimate together with  \eqref{Min12} clearly ensure that the family 
$\left\{\psi_{p,\mu},~0<\mu\ll 1\right\}$ 
is bounded in $H^{1}(\mathbb{R})$.

% CONTROLE DES MINIMISEURS DANS H^{2}
\textit{Step 4.} $H^{2}$-boundedness of the family $\{\psi_{p,\mu},~0<\mu< 1\}$ .

We will prove that the family $\{\psi_{p,\mu},~0<\mu< 1\}$ is bounded 
in $H^{2}(\mathbb{R})$. Applying the Fourier transformation
on equation \eqref{min010}, we obtain for all $\xi\in\mathbb{R}$,
\begin{equation}\label{Min14}
\left(1+\left|\xi\right|^{2}+\mu\left|\xi\right|^{4}\right)\widehat{\psi_{p,\mu}}(\xi)=
\frac{\alpha_{p,\mu}}{p+1}\widehat{\psi^{p+1}_{\mu}}(\xi).
\end{equation}
Since $\{\alpha_{p,\mu},~0<\mu< 1\}$ is bounded (from Step $3$), 
it suffices to prove that  
the family $\{\widehat{\psi^{p+1}_{p,\mu}},~0<\mu< 1\}$ 
is bounded in $L^{2}(\mathbb{R})$, 
and by the Plancherel-Parseval identity, 
this is equivalent to prove that the family
$\{\psi_{p,\mu},~0<\mu< 1\}$ is bounded in $L^{2p+2}(\mathbb{R})$. 
Thus, we compute
\begin{align}\label{Min15}
\left\|\psi_{p,\mu}\right\|^{2p+2}_{L^{2p+2}(\mathbb{R})}
&\le\left\|\psi_{p,\mu}\right\|^{2p}_{L^{\infty}(\mathbb{R})}
\left\|\psi_{p,\mu}\right\|^{2}_{L^{2}(\mathbb{R})}
\nonumber\\
&\le 2^{-p}\left\|\psi_{p,\mu}\right\|^{2p+2}_{H^{1}(\mathbb{R})},
\end{align}
where we use the Sobolev embedding of $H^{1}(\mathbb{R})$
into $C^{0}(\mathbb{R})$ (with the best constant).
Then, using that $\{\psi_{p,\mu},~0<\mu< 1\}$ 
is bounded in $H^{1}(\mathbb{R})$ (from Step $1$), and \eqref{Min15}, 
we obtain the desired result. Note that this boundedness, together with the fact that $ \varphi_{1,p,0} $ is a solution to the constraint minimization problem $ S^{\beta_p}_{p,0}$,  ensures that
\begin{equation}\label{Min44}
\liminf_{\mu\to 0} I_{\mu}(\psi_{p,\mu})=\liminf_{\mu\to 0} I_{0}(\psi_{p,\mu})\ge I_0(\varphi_{1,p,0}) \;.
\end{equation}
Gathering this last estimate with \eqref{Min12}-\eqref{Min13}, we get 
\begin{equation}\label{Min144}
 \quad \lim_{\mu\to 0} S^{\beta_p}_{p,\mu}=\lim_{\mu\to 0}  I_{\mu}(\psi_{p,\mu})=\lim_{\mu\to 0} I_{0}(\psi_{p,\mu})=I_0(\varphi_{1,p,0})  \quad \text{ and } \quad \lim_{\mu\to 0} \alpha_{p,\mu}= 1\;.
\end{equation}

% CONVERGENCE VERS LE SOLITON DE gKdV
\textit{Step 5.}  Strong convergence in $H^{1}(\mathbb{R})$  as $\mu\rightarrow 0$ and construction of the family of solitary waves $\{\varphi_{1,p,\mu},~0<\mu< 1\}$.

Let $(\mu_{n})_{n\ge 1}\subset\mathbb{R^{*}}_{+}$ be a sequence
which decreases to $0^{+}$, 
and for simplicity, let us rename $\psi_{p,\mu_{n}}=\psi_{n}$
and $\alpha_{p,\mu_{n}}=\alpha_{p,n}$. 
Then, there exist a sub-sequence, still denoted by $(\psi_{n})_{n\ge 1}$,
that converges weakly in $ H^2_e(\R) $ towards a function  $\psi_{\infty}\in H^2_e(\R)$. Passing to the limit in \eqref{min010} by making use of  \eqref{Min144}, 
we get that $  \psi_{\infty}$ satisfies \eqref{Int4} with $ c=1$.  
Note that, to pass to the limit on the nonlinear 
part of equation \eqref{min010} we deal as follows: using the pointwise 
convergence \eqref{min5-2-0}, the  $ H^2$-boundedness  
and the continuity of the function: $y\mapsto y^{p+1}$, 
one can check that 
\begin{equation}\label{min0004}
\psi^{p+1}_{n}\rightarrow\psi^{p+1}_{\infty}
~~\text{a.e. on}~~\mathbb{R}.
\end{equation}
Then, applying Lebesgue's Dominated Convergence Theorem, 
it holds for all $\phi\in H^{2}(\mathbb{R})$,
\begin{equation}\label{min0005}
\lim_{n\to+\infty}\langle\psi^{p+1}_{n},\phi\rangle_{L^{2}}
=\langle\psi^{p+1}_{\infty},\phi\rangle_{L^{2}}.
\end{equation}
Now, the uniqueness result (up to symmetries) of the solutions to \eqref{Int4} ensures that $ \psi_\infty=\varphi_{1,p,0} $ whenever $ p $ is even and 
$ \psi_\infty =\mp \varphi_{1,p,0}$ whenever $ p$ is odd. Since, according to \eqref{Min144},
$$
\lim_{n\to +\infty} I_0(\psi_n)= I_0(\mp \varphi_{1,p,0})=I_0(\psi_\infty) ,
$$
it follows that $ (\psi_n)_{n\ge 0} $ converges actually strongly to $\psi_\infty $ in $ H^1(\R) $ and thus $ K_p(\psi_\infty)=\beta_p>0 $. This ensures that in all cases, 
 $ \psi_\infty=\varphi_{1,p,0} $ and thus 
 \begin{equation}\label{gf}
 \psi_\mu \to \varphi_{1,n,0} \quad \text{ in } H^1(\R) \quad \text{ as } \mu\to 0 \; .
 \end{equation}

Finally, setting  
$$
\varphi_{1,p,\mu}=\alpha^{1/p}_{p,\mu}\psi_{p,\mu},\quad \forall n\in \N^* \; .
$$
 It is easy to check that for any $ p\in \N^* $ and any $ \mu>0 $, $\varphi_{1,p,\mu}$ is a solution of \eqref{Int3} with $ c=1 $. Moreover, according to
 Step 4., \eqref{Min144} and \eqref{gf}, the family  $\{\varphi_{1,p,\mu}, \quad 0<\mu < 1\} $ is bounded in $H^{2}(\mathbb{R})$ and
  \begin{equation}\label{gf2}
 \varphi_{1,p,\mu} \to \varphi_{1,p,0} \quad \text{ in } H^1(\R) \quad \text{ as } \mu\to 0 \; .
 \end{equation}
 This completes  the proof of the lemma.
\hfill $ \square $ \vspace*{2mm} 

% ALTERNATIVE
% More simply, using the strong convergence 
% in $H^{1}_{loc}(\mathbb{R})$ 
% (see \eqref{min5-2-1}), 
% and the uniform exponential decay properties of $(\psi_{n})_{n\ge 1}$ 
% and $\varphi_{c,p}$ (see \eqref{min0002}),
% we directly get the strong convergence in $H^{1}(\mathbb{R})$. 

% UNICITE DES SOLITONS PAIRS DE gKW POUR P ENTIER
\subsection{Uniqueness of even ground states}

 Note that the construction of the $\varphi_{1,p,\mu} $ ensures that they are even ground state solutions (in the sense of Definition \ref{groundstate}) to \eqref{Int4} with $c=1 $. 
 The following lemma ensures that for $ \mu$ small enough, the even ground state solution to  \eqref{Int4}, with $ c=1$,  is unique.
\begin{Lem}[Uniqueness of $\varphi_{1,p,\mu}$ for $0<\mu\ll 1$]\label{Uniq}
Let $p\in \{1,2,3\} $ be fixed and set $ c=1$. There exists $ \delta'_{p}>0$ such that for all $\mu \in ]0,\delta_p'[ $, there exists a unique even ground state solution
  $\varphi_{1,p,\mu}$ to equation \eqref{Int3}.
Moreover, the map $\mu\mapsto\varphi_{1,p,\mu}$ is-continuous from $ ]0,\delta_p'[ $ into $ H^1(\R) $.
\end{Lem}

% PREUVE
\textbf{Proof.} We fix $ p\in \{1,2,3\} $. 
We first claim that for any $ \varepsilon>0 $ there exists $ \mu_\varepsilon >0 $ such that for any $ 0<\mu <\mu_\varepsilon$, any even  ground state solution $ \varphi $ 
 to \eqref{Int3} with $ c=1$ satisfies
\begin{equation}\label{ds1}
\|\varphi-\varphi_{1,p,0} \|_{H^1} <\varepsilon \; .
\end{equation}
To prove this claim we proceed by contradiction, assuming that  there exist $ \varepsilon_0>0 $, a sequence $ (\mu_n)_{n\ge 0}$ of positive real numbers that converges to $ 0 $ and a sequence  $ (\varphi_n)_{n\ge 0} $ of ground state solutions to 
 \eqref{Int3} with $ c=1 $ and $ \mu=\mu_n$ such that 
 \begin{equation}\label{ds2}
 \|\varphi_n -\varphi_{1,p,0}\|_{H^1} \ge \varepsilon_0 \, , \quad \forall n\in \N\; .
 \end{equation}
 From \eqref{Int3} we infer that $ K_p(\varphi_n)=\frac{2}{p+2}I_{\mu_n}(\varphi_n)>0 $. We set 
 $$
 \psi_n=\Bigl(\frac{K_p(\varphi_{1,p,0})}{K_p(\varphi_n)}\Bigr)^{\frac{1}{p+2}} \varphi_n 
 $$
 so that $ K_p(\psi_n)=K_p(\varphi_{1,p,0})=\beta_p $.  $(\psi_n)_{n\ge 0} $ is thus a sequence of solutions to $ S^{\beta_p}_{p,\mu_n} $ and as in the proof of Lemma \ref{Exist-Conv} (Step 5.), it follows that $(\psi)_{n\ge 0} $ is bounded in $ H^2(\R) $ and $ \psi_n \to \varphi_{1,p,0} $ in $ H^1(\R) $. Therefore,
 $$
 \varphi_n=\Bigl( \frac{2I_{\mu}(\psi_{n})}{(p+2)\beta_{p}}\Bigr)^{1/p} \psi_n \to \varphi_{1,p,0} \quad \text{ in } H^1(\R) 
 $$
 which contradicts \eqref{ds2}.

Let   $\varphi_{1,p,\mu}$ and $\tilde{\varphi}_{1,p,\mu}$  be two ground
 ground states of \eqref{Int3}  with $ c=1 $ and $ \mu>0$. 
 Following the idea of Kenig et al. \cite{MR2859931} 
(Proposition $3$), we will prove, arguing by contradiction, 
that $\tilde{w}=\varphi_{1,p,\mu}-\tilde{\varphi}_{1,p,\mu}=0$  as soon as $ \mu$ is small enough. Note that, on account of \eqref{ds1}, it holds 
$\|\tilde{w}\|_{H^1}\le \varepsilon(\mu) $ with $ \varepsilon(y) \to 0 $ as $ y \to 0 $.
We set $ F(x)=\frac{1}{p+1}x^{p+1} $ so that $ \tilde{w} $ satisfies 

% PREUVE UNICITE DE SOLITONS PAIRS DE gKW POUR P ENTIER
\begin{equation}\label{Min019}
\mu\partial^{4}_{x}\tilde{w}-\partial^{2}_{x} \tilde{w}
+ \tilde{w} =F(\varphi_{1,p,\mu})-F(\tilde{\varphi}_{1,p,\mu})
\end{equation}
Denoting by   $\mathcal{L}_{1,p,\mu}$ the operator defined in Proposition \ref{Anne} with   $\phi_{c,p,\mu}=\varphi_{1,p,\mu}$, it holds
$$
\mathcal{L}_{1,p,\mu}= \mu \partial_x^4 -\partial_x^2 +1 -F'(\varphi_{1,p,\mu}) 
$$
and thus 
$$
\mathcal{L}_{1,p,\mu} \tilde{w}= F(\varphi_{1,p,\mu})-F(\tilde{\varphi}_{1,p,\mu})+F'(\varphi_{1,p,\mu}) (\varphi_{1,p,\mu}-\tilde{\varphi}_{1,p,\mu})
=\frac{F''(W_{p,\mu})}{2} \tilde{w}^2 \, , 
$$
with, for all $x\in \R $, $W_{p,\mu}(x)\in [\varphi_{1,p,\mu}(x), \tilde{\varphi}_{1,p,\mu}(x)] $. Therefore, assuming that $ \tilde{w} \neq 0 $ and setting $ w=\frac{\tilde{w}}{\|\tilde{w}\|_{H^1}} $,
 we get 
\begin{equation}\label{Min023}
\left\|\mathcal{L}_{1,p,\mu}w\right\|_{L^{2}(\mathbb{R})}
\lesssim
\frac{1}{\|\tilde{w}\|_{H^1}}\|\tilde{w}\|_{L^2}\|\tilde{w}\|_{H^1}\lesssim  \|\tilde{w}\|_{L^2} \|w\|_{H^1} \, 
\end{equation}
so that
\begin{align}\label{Min025}
\left|\left\langle\mathcal{L}_{1,p,0}w,w\right\rangle_{L^{2}}\right|
&=\left|\left\langle \mathcal{L}_{1,p,\mu}w,w\right\rangle_{L^{2}}+\left\langle w^2,\varphi_{1,p,\mu}^p -\varphi_{1,p,0}^p\right\rangle_{L^{2}}
-\mu \left\|\partial^{2}_{x}w\right\|^{2}_{L^{2}(\mathbb{R})}\right|
\nonumber\\
&\lesssim  \varepsilon(\mu)+\mu \; .
\end{align}
On the other hand, it is well-known that for $ p\in \{1,2,3\}$, $ \frac{d}{dc}_{\vert_{c=1}} \|\varphi_{c,p,0}\|_{L^2}^2 <0 $ which ensures that there exists $ \alpha_p>0 $ such that 
$$
\langle {\mathcal L}_{1,p,0} \phi,\phi\rangle \ge \alpha_p \|\phi\|_{H^1}^2
$$
for all $ \phi\in H^1(\R) $ that satisfies the orthogonality conditions : $\langle \phi, \varphi_{1,p,0}'\rangle= \langle \phi, \varphi_{1,p,0}\rangle=0 $.
\noindent
Now, since $w $ is even and $\varphi'_{1,p,\mu}$ is odd,
we have the first orthogonality condition:
\begin{equation}\label{Min026}
\left\langle w_{n},\varphi'_{c,p}\right\rangle_{L^{2}}=0.
\end{equation}
Moreover, using that $ {\mathcal L}_{1,p,0} {\partial_{c}}_{\vert c=1} \varphi_{c,p,0}=\varphi_{1,p,0}$ we get 
\begin{align}\label{Min0027}
\left|\left\langle w,\varphi_{1,p,0}\right\rangle_{L^{2}}\right|
&=\left|\left\langle w , {\mathcal L}_{1,p,0} {\partial_{c}}_{\vert c=1} \varphi_{c,p,0}
\right\rangle_{L^{2}}\right|
\nonumber\\
&=\left|\left\langle {\mathcal L}_{1,p,0}  w , {\partial_{c}}_{\vert c=1} \varphi_{c,p,0}
\right\rangle_{L^{2}}\right|
\nonumber\\
&\le\left|\left\langle\mathcal{L}_{1,p,\mu}w, {\partial_{c}}_{\vert c=1} \varphi_{c,p,0}
\right\rangle_{L^{2}}\right|
+\left|\left\langle w  {\partial_{c}}_{\vert c=1} \varphi_{c,p,0},\varphi^{p}_{c,p,\mu}
-\varphi^{p}_{1,p,0}\right\rangle_{L^{2}}
\right|
\nonumber\\
&\hspace{1cm}+\mu \left|\left\langle
\partial^{2}_{x}w ,\partial^{2}_{x} {\partial_{c}}_{\vert c=1} \varphi_{c,p,0}
\right\rangle_{L^{2}}\right|
\nonumber\\
&\lesssim  (\varepsilon(\mu)+\mu) \|w\|_{H^1}
\end{align}
so that  for $ \mu$ small enough, $ w$ is almost orthogonal to $ \varphi_{1,p,0}$. Therefore, proceeding as in Subsection \ref{Stab-gKW},
we infer that for $ \mu >0 $ small enough, 
$$
\langle {\mathcal L}_{1,p,0} w,w\rangle \ge \frac{\alpha_p}{2} \|w\|_{H^1}^2= \frac{\alpha_p}{2},
$$
which contradicts \eqref{Min025} and proves the uniqueness result for $ 0<\mu<\delta_p' $. Finally, once the uniqueness is proved, the convergence of 
 $ \varphi_{1,p,\mu} $ towards $  \varphi_{1,p,\mu_0} $ as $ \mu\to \mu_0 $ can be proved exactly as the convergence of  $ \varphi_{1,p,\mu} $ towards $  \varphi_{1,p,0} $ as $ \mu\to 0 $.
\hfill $ \square $ \vspace*{2mm} \\
Now, we notice that $ v $ satisfies \eqref{Int3} with $ p\in \N^*$, $c=1 $ and $ \mu>0 $  if and only $ v_\mu=\mu^{1/p}  v(\sqrt{\mu}\,\, \cdot) $ satisfies
  \eqref{Int3} with $ p\in \N^*$, $c=\mu $ and $ \mu=1 $. The following 
 uniqueness result on the ground states to \eqref{Int3} with $ \mu=1$ is thus a direct consequence of Lemma \ref{Uniq}.

% STABILITE DES SOLITONS LENTS DE gKW
\subsection{End of the proof of Theorem \ref{Main-Result}}
Let $ p\in \{1,2,3\} $. 
It remains to prove the orbital stability result. Since $u $ is a solution to \eqref{Int1}  if and only if $ u_\mu(t,x)= \mu^{1/p} u(\mu^{3/2} t, \sqrt{\mu}x)
$ is a solution to \eqref{Int1} with $ \mu=1 $, the orbital stability of $ \varphi_{\mu,p,1} $ is equivalent to the orbital stability of $ \varphi_{c,p,\mu} $. To prove the orbital stability of this last ground state for $ \mu>0 $ small enough, we rely on a continuity argument as in the preceding section.
 First, by the $H^1$-convergence result \eqref{Min04}, it is easy to check that for $ \mu >0 $ small enough, any $ v\in H^2(\R) $ satisfying 
 $ \langle v, \varphi_{1,p,\mu}\rangle= \langle v, \varphi_{1,p,\mu}'\rangle= 0 $ is almost orthogonal in $ L^2(\R) $ to $ \varphi_{1,p,0} $ and 
  $ \varphi_{1,p,0}'$. Second, recalling that, for $ p\in \{1,2,3\} $, $ \frac{d}{dc}_{\vert_{c=1}} \|\varphi_{c,p,0}\|_{L^2}^2 <0 $
 and arguing exactly as in Subsection
\ref{Stab-gKW}, we obtain that $\mathcal{L}_{1,p,0}$
is coercive in $H^{1}(\mathbb{R})$ under these almost orthogonality conditions. It follows that
\begin{align}\label{Main3}
\left\langle \mathcal{L}_{1,p,\mu}v,v\right\rangle_{L^{2}}&=
\mu\left\|\partial^{2}_{x}v\right\|^{2}_{L^{2}(\mathbb{R})}
+\left\langle\mathcal{L}_{1,p}v,v\right\rangle_{L^{2}}
-\left\langle v^{2},\varphi^{p}_{1,p,\mu}-\varphi^{p}_{1,p,0}
\right\rangle_{L^{2}}
\nonumber\\
&\ge\mu\left\|\partial^{2}_{x}v\right\|^{2}_{L^{2}(\mathbb{R})}
+\frac{\alpha_{p}}{2}\left\|v\right\|^{2}_{H^{1}(\mathbb{R})}
-\varepsilon(\mu) \left\|v\right\|^{2}_{L^{2}(\mathbb{R})}
\nonumber\\
&\ge\min\{\mu,\frac{\alpha_p}{2}-\varepsilon(\mu)\}
\left\|v\right\|^{2}_{H^{2}(\mathbb{R})},
\end{align}
with  $\alpha_{p}>0$ and $ \varepsilon(y) \to 0 $ as $ y\to 0 $. 
By applying Proposition \ref{Anne}, 
we  thus get that  there exists $ \delta_p>0 $ such that for all $ \mu \in]0,\delta_p[$,  $ \varphi_{1,p,\mu} $ is orbitally stable. This completes the proof of the theorem.
\hfill $ \square $ \vspace*{2mm}

% REFERENCES 
\newpage
\bibliographystyle{plain}
\bibliography{mabiblio}

\begin{thebibliography}{10}

\bibitem{MR1044731}
L.~Abdelouhab, J.~L. Bona, M.~Felland, and J.-C. Saut.
\newblock Nonlocal models for nonlinear, dispersive waves.
\newblock {\em Phys. D}, 40(3):360--392, 1989.

\bibitem{MR1151253}
J.~P. Albert.
\newblock Positivity properties and stability of solitary-wave solutions of
  model equations for long waves.
\newblock {\em Comm. Partial Differential Equations}, 17(1-2):1--22, 1992.

\bibitem{MR0338584}
T.~B. Benjamin.
\newblock The stability of solitary waves.
\newblock {\em Proc. Roy. Soc. (London) Ser. A}, 328:153--183, 1972.

\bibitem{MR0386438}
J.~Bona.
\newblock On the stability theory of solitary waves.
\newblock {\em Proc. Roy. Soc. London Ser. A}, 344(1638):363--374, 1975.

\bibitem{MR1946769}
Thomas~J. Bridges, Gianne Derks, and Georg Gottwald.
\newblock Stability and instability of solitary waves of the fifth-order
  {K}d{V} equation: a numerical framework.
\newblock {\em Phys. D}, 172(1-4):190--216, 2002.

\bibitem{MR3461359}
R.~C{\^o}te, C.~Mu{\~n}oz, D.~Pilod, and G.~Simpson.
\newblock Asymptotic {S}tability of {H}igh-dimensional {Z}akharov--{K}uznetsov
  {S}olitons.
\newblock {\em Arch. Ration. Mech. Anal.}, 220(2):639--710, 2016.

\bibitem{Anne}
A.~de~Bouard.
\newblock Equations dispersives non lin\'eaires.
\newblock {\em Cours \`a l'Ecole d'Et\'e de Math\'ematiques. Institut Fourier,
  Grenoble.}, pages 1--35, 2005.

\bibitem{MR1427885}
B.~Dey, A.~Khare, and C.~N. Kumar.
\newblock Stationary solitons of the fifth order {K}d{V}-type. {E}quations and
  their stabilization.
\newblock {\em Phys. Lett. A}, 223(6):449--452, 1996.

\bibitem{MR1372681}
V.~I. Karpman.
\newblock Stabilization of soliton instabilities by higher order dispersion:
  {K}d{V}-type equations.
\newblock {\em Phys. Lett. A}, 210(1-2):77--84, 1996.

\bibitem{MR1837757}
V.~I. Karpman.
\newblock Stability of solitons described by {K}d{V}-type equations with
  higher-order dispersion and nonlinearity (analytical results).
\newblock {\em Phys. Lett. A}, 284(6):238--246, 2001.

\bibitem{Kawa}
T.~Kawahara.
\newblock Oscillatory solitary waves in dispersive media.
\newblock {\em Phys. Soc. Japan}, 33:260--264, 1972.

\bibitem{MR2859931}
C.~E. Kenig, Y.~Martel, and L.~Robbiano.
\newblock Local well-posedness and blow-up in the energy space for a class of
  {$L^2$} critical dispersion generalized {B}enjamin-{O}no equations.
\newblock {\em Ann. Inst. H. Poincar\'e Anal. Non Lin\'eaire}, 28(6):853--887,
  2011.

\bibitem{MR1667522}
S.~P. Levandosky.
\newblock A stability analysis of fifth-order water wave models.
\newblock {\em Phys. D}, 125(3-4):222--240, 1999.

\bibitem{MR2332504}
S.~P. Levandosky.
\newblock Stability of solitary waves of a fifth-order water wave model.
\newblock {\em Phys. D}, 227(2):162--172, 2007.

\bibitem{MR778970}
P.-L. Lions.
\newblock The concentration-compactness principle in the calculus of
  variations. {T}he locally compact case. {I}.
\newblock {\em Ann. Inst. H. Poincar\'e Anal. Non Lin\'eaire}, 1(2):109--145,
  1984.

\bibitem{MR0232968}
W.~Magnus, F.~Oberhettinger, and R.~Soni.
\newblock {\em Formulas and theorems for the special functions of mathematical
  physics}.
\newblock Third enlarged edition. Die Grundlehren der mathematischen
  Wissenschaften, Band 52. Springer-Verlag New York, Inc., New York, 1966.

\bibitem{MR2215276}
M.~Mari{\c{s}}.
\newblock Global branches of travelling-waves to a
  {G}ross-{P}itaevskii-{S}chr\"odinger system in one dimension.
\newblock {\em SIAM J. Math. Anal.}, 37(5):1535--1559 (electronic), 2006.

\bibitem{MR1958041}
Angulo P.J.
\newblock On the instability of solitary-wave solutions for fifth-order water
  wave models.
\newblock {\em Electron. J. Differential Equations}, pages No. 6, 18 pp.
  (electronic), 2003.

\bibitem{MR783974}
M.~Weinstein.
\newblock Modulational stability of ground states of nonlinear {S}chr\"odinger
  equations.
\newblock {\em SIAM J. Math. Anal.}, 16(3):472--491, 1985.

\end{thebibliography}

% FIN DU DOCUMENT
\end{document}